\documentclass[11pt,letterpaper,reqno]{amsart}
\usepackage{epsfig}

\def\be{\begin{equation}}
\def\ee{\end{equation}}
\def\bea{\begin{eqnarray}}
\def\eea{\end{eqnarray}}
\def\bes{\begin{eqnarray*}}
\def\ees{\end{eqnarray*}}

\def\nn{\nonumber}
\def\lb{\label}
\def\bs{\setminus}
\def\pt{\partial}

\def\T{{\mathcal{T}}}
\def\H{{\mathcal{H}}}
\def\R{{\bf R}}
\def\C{{\bf C}}
\def\Z{{\bf Z}}

\def\N{{\bf N}}
\def\U{{\bf U}}

\def\Q{{\bf Q}}
\def\T{{\bf T}}

\def\Sg{{\Sigma}}

\def\aa{{\alpha}}
\def\bb{{\beta}}
\def\ga{{\gamma}}

\def\th{{\theta}}

\def\om{{\omega}}
\def\Om{{\Omega}}
\def\ep{{\epsilon}}
\def\lm{{\lambda}}
\def\Lm{{\Lambda}}
\def\dl{{\delta}}

\def\sg{{\sigma}}
\def\dm{{\diamond}}
\def\Sg{{\Sigma}}
\def\vf{{\varphi}}

\def\<{{\langle}}
\def\>{{\rangle}}
\def\T{{\mathcal{T}}}

\def\Nn{{\mathcal{N}}}

\def\mul{{\rm mul}}

\def\crit{{\rm crit}}

\def\Sp{{\rm Sp}}

\def\hb{\vrule height0.18cm width0.14cm $\,$}

\title[Generalized common index jump theorem with applications]{Generalized common index jump theorem with applications to closed characteristics on star-shaped hypersurfaces and beyond}

\author[Huagui Duan]{Huagui Duan$^1$}
\thanks{$^1$ Partially supported by National Key R\&D Program of China (Grant No. 2020YFA0713300)
and NSFC (Nos. 12271268, 11671215 and 11790271), LPMC of MOE of China and Nankai University.}
\address{(Huagui Duan) School of Mathematical Sciences and LPMC, Nankai University, Tianjin 300071, The People's Republic of China}
\email{duanhg@nankai.edu.cn.}

\author[Hui Liu]{Hui Liu$^2$}
\thanks{$^2$ Partially supported by NSFC (Nos. 12022111 and 11771341).}
\address{(Hui Liu) School of Mathematics and Statistics, Wuhan University,
Wuhan 430072, Hubei, The People's Republic of China}
\email{huiliu00031514@whu.edu.cn.}

\author[Yiming Long]{Yiming Long$^3$}
\thanks{$^3$ Partially supported by National Key R$\&$D Program of China (Grant No. 2020YFA0713300),
NSFC (Nos. 11131004, 11671215 and 11790271), LPMC of Ministry of Education of China, Nankai University, Wenzhong
Foundation, and Nankai Zhide Foundation.}
\address{(Yiming Long) Chern Institute of Mathematics and LPMC, Nankai University, Tianjin 300071, The People's Republic of China}
\email{longym@nankai.edu.cn.}

\author[Wei Wang]{Wei Wang$^4$}
\thanks{$^4$ Partially supported by NSFC (Nos. 12025101 and 11431001).}
\address{(Wei Wang) Key Laboratory of Pure and Applied Mathematics,
School of Mathematical Sciences, Peking University, Beijing 100871,
The People's Republic of China}
\email{wangwei@math.pku.edu.cn.}

\begin{document}
\maketitle

\begin{abstract}
{\it In this paper, we first generalize the common index jump theorem of Long-Zhu in 2002 and Duan-Long-Wang
in 2016 to the case where the mean indices of symplectic paths are not required to be all positive. As applications,
we study closed characteristics on compact star-shaped hypersurfaces in $\R^{2n}$, when both positive and negative mean
indices may appear simultaneously. Specially we establish the existence of at least $n$ geometrically distinct closed characteristics on every
compact non-degenerate perfect star-shaped hypersurface $\Sg$ in $\R^{2n}$ provided that every prime closed
characteristic possesses nonzero mean index. Furthermore, in the case of $\R^6$ we remove the nonzero mean index
condition by showing that the existence of only finitely many geometrically distinct closed characteristics implies that each of them must
possess nonzero mean index. We also generalize the above results about closed characteristics on non-degenerate star-shaped hypersurfaces to
closed Reeb orbits of non-degenerate contact forms on a broad class of prequantization
bundles. }
\end{abstract}



\renewcommand{\theequation}{\thesection.\arabic{equation}}
\renewcommand{\thefigure}{\thesection.\arabic{figure}}

\setcounter{figure}{0}
\setcounter{equation}{0}
\section{Introduction and main results}

There are three goals in this paper. The first one is to generalize the common index jump theorem (CIJT for short
below). This theorem is the first result which exhibits certain common index property of iterates of more than one
but finitely many symplectic matrix paths. This theorem was proved by Y. Long and C. Zhu in 2002
(Theorems 4.1 and 4.3 of \cite{LoZ}, cf. also Theorems 11.1.1 and 11.2.1 of \cite{Lon4}).
In the last twenty years, this theorem has been applied to study various problems, provides an important method in
the studies of the multiplicity and stability of closed characteristics in symplectic and contact dynamics. Specially
it is one of the few methods which work effectively for such studies when the dimension of the symplectic manifold is
not only $4$. In \cite{Wan2} and \cite{Wan3} of 2016, W. Wang found certain useful symmetric property in CIJT.
In 2016 also, an important extension of this CIJT was established by H. Duan, Y. Long and W. Wang (Theorem 3.5 of \cite{DLW}),
which gives an enhanced version of the CIJT by giving precise indices for iterates of related symplectic paths near the
carefully chosen iterates in the study. This enhanced CIJT has been used to establish sharp estimates on the
multiplicity and stability of prime closed geodesics on compact simply-connected bumpy Finsler manifolds whose loop
spaces possess bounded Betti number sequences, provided the number of prime closed geodesics is finite and each of them
possesses non-zero Morse index, which implies the positivity of their mean indices in \cite{DLW}. This enhanced CIJT
has also been applied to the studies of closed characteristics on star-shaped hypersurfaces in $\R^{2n}$ by H. Duan,
H. Liu, Y. Long and W. Wang in \cite{DLLW} and on other contact manifolds by V. Ginzburg, B. Z. G\"urel and L. Macarini
in \cite{GGM} for example, under the assumption that all the prime closed characteristics possess positive (or negative)
mean indices. Note that in Theorems 5.1 and 5.2 as well as Corollaries 5.3 and 5.4 of \cite{GGu}, V. Ginzburg and B. Z.
G\"urel extended the enhanced CIJT of \cite{DLW} to the case admitting the mean indices of all symplectic paths being
nonzero via a new index recurrence arguments. Subsequently, based on \cite{GGu}, V. Ginzburg, B. Z. G\"urel and L.
Macarini in Theorem 4.1 of \cite{GGM} further studied the enhanced CIJT of \cite{DLW} for strongly non-degenerate
symplectic paths with positive mean indices. Note that based on the theorems these authors gave interesting results about
Reeb orbits on contact manifolds. But their extensions of the enhanced CIJT missed the precise values of indices of some
crucial iterates as those listed in Theorem 3.6 below, and missed the symmetric property of the CIJT produced by the
vertices in the cube $[0,1]^l$ in the proof of CIJT discovered first in \cite{LoZ} and then in \cite{DLW} as those in the
Step 2 of the proof of Theorem 1.2 below where the two opposite vertices are used. Note that such missing might be due
to the index recurrence arguments. Note also that these missing properties are very crucial in our study in the current
paper and for the future study, more precisely in order to get sharp estimates on multiplicities of periodic orbits, we do
need to compute Morse type number quantities accurately and to apply the mentioned symmetric property of CIJM. Thus the
first goal of this paper is to further extend this enhanced CIJT of \cite{DLW} to the case that there exist prime
symplectic paths possessing positive as well as negative mean indices simultaneously and the above mentioned information
can be reserved at the same time by generalizing the method of \cite{DLW}.

The main idea in the proof of the CIJT in \cite{LoZ} is to show that there exist large suitable iterate of each
path among the given finitely many prime symplectic paths such that the corresponding enlarged index intervals
at their these iterates possess a non-empty common intersection interval which is sufficiently large to contain
certain integers, and the number of these integers yields a lower bound estimate for the number of these prime
symplectic paths, provided all of these paths possess positive mean indices. Note that when prime closed characteristics
on a compact star-shaped hypersurface in $\R^{2n}$ are considered, in general some of them may possess positive mean
indices and the others possess negative mean indices. Then the iterated index sequences of them may tend to positive
infinity as well as negative infinity simultaneously, and consequently it seems impossible to apply the CIJT to get a
common intersection interval of the iterated enlarged index intervals of all the prime symplectic paths. To overcome
this difficulty, suggested by the resonance identity for a star-shaped hypersurface in $\R^{2n}$ with only finitely
many prime closed characteristics which was proved by H. Liu, Y. Long and W. Wang in \cite{LLW} of 2014 (cf. Theorem 2.3
below), we can consider iterates of all these prime symplectic paths together by adding some $-1$ parameter to the
negative mean indices, and we have improved the common index jump theorem (Theorems 4.1 and 4.3 of \cite{LoZ}
as well as Theorems 11.1.1 and 11.2.1 of \cite{Lon4}) to the new enhanced common index jump Theorems 3.4 and 3.6 below to
deal with mixed mean indices. This improvement allows our Theorem 1.2 below to deal with closed characteristics on
non-degenerate star-shaped hypersurfaces when their mean indices are non-zero.

The second goal of this paper is to apply this generalized enhanced common index jump Theorem to study the multiplicity
and stability of closed characteristics on compact star-shaped hypersurfaces in $\R^{2n}$. In this paper, we let
$\Sigma$ be always a $C^3$ compact hypersurface in $\R^{2n}$ strictly star-shaped with respect to the origin, i.e.,
the tangent hyperplane at any $x\in\Sigma$ does not intersect the origin. We denote the set of all such hypersurfaces
by $\H_{st}(2n)$, and denote by $\H_{con}(2n)$ the subset of $\H_{st}(2n)$ which consists of all strictly convex
hypersurfaces. We consider closed characteristics $(\tau, y)$ on $\Sigma$, which are solutions of the following
problem
\be
\left\{\begin{array}{ll}\dot{y}=JN_\Sigma(y), \\
               y(\tau)=y(0),  \lb{1.1}
\end{array}\right.
\ee
where $J=\left(\begin{array}{ll}0 &-I_n\\
        I_n  & 0\end{array}\right)$, $I_n$ is the identity matrix in $\R^n$, $\tau>0$, $N_\Sigma(y)$ is the outward
normal vector of $\Sigma$ at $y$ normalized by the condition $N_\Sigma(y)\cdot y=1$. Here $a\cdot b$ denotes
the standard inner product of $a, b\in\R^{2n}$. A closed characteristic $(\tau, y)$ is {\it prime}, if $\tau$
is the minimal period of $y$. Two closed characteristics $(\tau, y)$ and $(\sigma, z)$ are {\it geometrically
distinct}, if $y(\R)\not= z(\R)$. We denote by $\T(\Sigma)$ the set of geometrically distinct
closed characteristics $(\tau, y)$ on $\Sigma\in\mathcal{H}_{st}(2n)$. A closed characteristic
$(\tau,y)$ is {\it non-degenerate} if $1$ is a Floquet multiplier of $y$ of precisely algebraic multiplicity
$2$; {\it hyperbolic} if $1$ is a double Floquet multiplier of it and all the other Floquet multipliers
are not on ${\bf U}=\{z\in {\bf C}\mid |z|=1\}$, i.e., the unit circle in the complex plane; {\it elliptic}
if all the Floquet multipliers of $y$ are on ${\bf U}$. We call a $\Sigma\in \mathcal{H}_{st}(2n)$ {\it
non-degenerate} if all the closed characteristics on $\Sigma$, together with all of their iterations, are
non-degenerate.

The study on closed characteristics in the global sense started in 1978, when the existence of at least one
closed characteristic was first established on any $\Sg\in\H_{st}(2n)$ by P. Rabinowitz in \cite{Rab1}
and on any $\Sg\in\H_{con}(2n)$ by A. Weinstein in \cite{Wei1} independently. Since then the existence of
multiple closed characteristics on $\Sg\in\H_{con}(2n)$ has been deeply studied by many mathematicians, for
example, studies in \cite{EkL1}, \cite{EkH1}, \cite{HWZ1}, \cite{Szu1}, \cite{LoZ}, \cite{LLZ}, \cite{Wan2},
\cite{Wan3}, \cite{WHL} as well as \cite{Lon4} and references therein.

For the star-shaped hypersurfaces, We are only aware of a few papers about the multiplicity of closed
characteristics. In \cite{Gir1} of 1984 and \cite{BLMR} of 1985, $\;^{\#}\T(\Sg)\ge n$ for $\Sg\in\H_{st}(2n)$
was proved under some pinching conditions. In \cite{Vit2} of 1989, C. Viterbo proved a generic existence
result for infinitely many closed characteristics on star-shaped hypersurfaces. In \cite{HuL} of 2002, X. Hu
and Y. Long proved that $\;^{\#}\T(\Sg)\ge 2$ for non-degenerate $\Sg\in \H_{st}(2n)$. In \cite{HWZ2} of 2003,
H. Hofer, K. Wysocki, and E. Zehnder proved that $\,^{\#}\T(\Sg)=2$ or $\infty$ holds for every non-degenerate
$\Sg\in\H_{st}(4)$ provided that all stable and unstable manifolds of the hyperbolic closed characteristics
on $\Sg$ intersect transversally. Furthermore, recently this alternative result was proved to be true for every
non-denernerate $\Sg\in\H_{st}(4)$ without the above transversal condition by D. Cristofaro-Gardiner, M.
Hutchings and D. Pomerleano in \cite{CGHP}. In \cite{CGH1} of 2016, D. Cristofaro-Gardiner and M. Hutchings
proved that $\;^{\#}\T(\Sg)\ge 2$ for every contact manifold $\Sg$ of dimension three. Later various proofs of
this result for star-shaped hypersurfaces have been given in \cite{GHHM}, \cite{LLo1} and \cite{GiG1}. There
are also some multiplicity results for closed orbits of dynamically convex Reeb flows, cf., \cite{GuK} and
\cite{AbM}.

On the stability problem, we refer the readers to \cite{Eke1}, \cite{DDE1}, \cite{Lon1}-\cite{Lon3}, \cite{LoZ},
\cite{WHL}, \cite{Wan1} for convex hypersurfaces and \cite{LiL}, \cite{LLo2}, \cite{CGHHL} for star-shaped
hypersurfaces. The following index perfect condition was first introduced by H. Duan, H. Liu, Y. Long and W. Wang
in \cite{DLLW} of 2018, for the star-shaped hypersurfaces in $\R^{2n}$ which is much weaker than the dynamically
convexity condition introduced by H. Hofer, K. Wysocki, and E. Zehnder in \cite{HWZ1} of 1998 (cf. also \cite{HWZ2}
in 2003).

\medskip

{\bf Definition 1.1.} {\it We call a compact star-shaped hypersurface $\Sigma$ in $\R^{2n}$ {\rm perfect}, if
for every prime closed characteristic $(\tau,y)$ on $\Sigma$, the Maslov-type index of each {\rm good} $m$-th
iterate $(m\tau,y)$ of $(\tau,y)$ with some $m\in\N$ satisfies $i(y,m)\not= -1$ if $n$ is even, or
$i(y,m)\not\in \{-2,-1,0\}$ if $n$ is odd. }

\medskip

Here an iterate $(m\tau, y)$ of a prime closed characteristic $(\tau,y)$ on $\Sigma$ with $m\in\N$ is called
{\it good}, if its Maslov-type index has the same parity as that of $(\tau,y)$, otherwise it is called {\it bad}.
Note that for a bad closed characteristic $(m\tau,y)$, the element $x^m$ corresponding to $(m\tau,y)$ satisfies
$\bb(x^m)=-1$ in Lemma 2.2 below, and consequently the equivariant critical module of the functional $F_{a,K}$
at $S^1\cdot x^m$ must be trivial, i.e., $x^m$ is homologically invisible and thus can be ignored in the Morse
theory study. This property was used first in Definition 4.8 and Remark 4.9 via Proposition 4.2 of \cite{LLW} to
compute the Euler characteristic $\hat{\chi}(y)$. Then this property was used in Section 3 of \cite{GGM} to
compute the local equivariant symplectic homology. (cf. also the condition (A) below Theorem 1.2). Here in the
current paper, we shall use this property in the computations of the Morse type numbers in Section 4 below.

In \cite{DLLW}, the authors proved specially the existence of at least $n$ closed characteristics on every
non-degenerate perfect star-shaped hypersurfaces in $\R^{2n}$ when every closed characteristic possesses positive
mean index. Then V. Ginzburg, B.Z. G\"{u}rel and L. Macarini in \cite{GGM} obtained the same multiplicity result
of closed Reeb orbits on contact manifolds under non-degenerate condition and the index perfect condition
introduced in \cite{DLLW} (i.e., the perfect condition given in Definition 1.1), provided the contact form $\aa$
is index-positive (or index-negative), i.e., all contractible periodic orbit $\gamma$ of
$\aa$ possess positive (or negative) mean index. Most recently V. Ginzburg and L. Macarini in \cite{GM} obtained
some optimal multiplicity results of closed Reeb orbits on symmetric contact spheres under the so called strong
dynamical convexity which extended results of C. Liu, Y. Long and C. Zhu in \cite{LLZ} of 2002.

When we consider prime closed characteristics on a compact star-shaped hypersurface in $\R^{2n}$, a priori the
mean indices of some of them can be non-positive. By the resonance identity for a star-shaped hypersurface in
$\R^{2n}$ with only finitely many prime closed characteristics proved by H. Liu, Y. Long and W. Wang in
\cite{LLW} of 2014 (cf. Theorem 2.3 below), at least one of these prime closed characteristics must possess
positive mean index, but some of the others may possess zero or negative mean indices. Thus even if we assume
that every prime closed characteristic possesses nonzero mean index, their iterated index sequences may still tend
to positive infinity as well as negative infinity simultaneously as we mentioned before. To overcome this difficulty,
by Theorem 2.3 below, we notice that the existence of some prime closed characteristics possessing positive mean
index can be crucially used to construct actually the mentioned common index intersection interval, and that the
prime closed characteristics with negative mean indices make no contributions to it and thus can be ignored in some
sense. This understanding is rigorously realized by applying our new generalized enhanced common index jump Theorem
3.6 together with the mentioned existence of prime closed characteristics possessing positive mean indices so that
we can deal with positive and negative mean indices simultaneously, and establish the following more general
existence result, provided the mean indices of all the prime closed characteristics are non-zero. Note that another key
observation in our proof is that the Morse inequalities still hold under non-zero mean index condition when the degrees
of Morse-type numbers and Betti numbers are carefully chosen (cf. (\ref{2.16}) and (\ref{2.17}) below).

\medskip

{\bf Theorem 1.2.} {\it Let $\Sigma$ be a compact non-degenerate perfect star-shaped hypersurface in $\R^{2n}$.
If every prime closed characteristic on $\Sigma$ possesses nonzero mean index, then there exist at least $n$
geometrically distinct closed characteristics. Furthermore, if the total number of prime closed characteristics on
$\Sigma$ is finite, then $\Sigma$ carries at least $n$ non-hyperbolic closed characteristics with even Maslov-type
indices when $n$ is even, and at least $n$ closed characteristics with odd Maslov-type indices and at least $(n-1)$
of them are non-hyperbolic when $n$ is odd.}

\medskip

Here we briefly explain the main ideas of the proof of Theorem 1.2 in three steps when $n$ is even, see Section 4 for details.

More precisely, in order to prove Theorem 1.2, we assume that there exist only finitely many prime closed characteristics $\{y_k\}_{k=1}^q$ with $\hat{i}(y_k)\neq 0$ on a non-degenerate perfect star-shaped hypersurface $\Sigma$ in $\R^{2n}$.

{\bf Step A.} Firstly, we apply the generalized common index jump theorem, i.e., Theorem 3.6, to get a $(q+1)$-tuple
$(N, m_1, \cdots, m_q) \in \N^{q+1}$ such that the Maslov-type indices $i(y_k^h)$ of $h=2m_k\pm m$ and $h=2m_k$-th
iterates of each prime closed characteristics $y_k$ satisfy the inequalities (\ref{4.10})-(\ref{4.14'}) with $1\le m\le \bar{m}$ for some suitably chosen $\bar{m}\in\N$. Here the perfect condition and $\hat{i}(y_k)\neq0$ with $1\le k\le q$ are used.

Then using these information on the these Maslov-type indices, we can compute
the alternating sum of Morse type numbers between $-2N-n-1$ and $2N-n-1$, i.e., $\sum_{p=-2N-n-1}^{2N-n-1}(-1)^p M_p$, which can be proved to be the
sum of $\sum_{k=1}^q2m_k\hat{\chi}_{y_k}=N$ and $(N_+^o - N_+^e)$, where $\hat{\chi}$ is the average Euler characteristic in Theorem 2.3 and $N_+^o,  N_+^e$ are defined by (\ref{4.21})-(\ref{4.22}).

Now by direct computations and the Morse inequality, we obtain
$$ N+ N_+^o - N_+^e = \sum_{p=-2N-n-1}^{2N-n-1}(-1)^p M_p \le \sum_{p=-2N-n-1}^{2N-n-1}(-1)^pb_p(\Lm M)
    = N - \frac{n}{2},  $$
i.e.,
$$ N_+^e \ge \frac{n}{2}. $$

{\bf Step B.} Note that the $(q+1)$-tuple $(N, m_1, \cdots, m_q)$ in Step A is chosen according to a vertex
$\chi$ of the cub $[0,1]^l$ given in Remark 3.7. Then for the vertex $\hat{\chi}=\hat{1}-\chi$ opposite
to $\chi$ in $[0,1]^l$ where $\hat{1}=(1,\cdots,1)$, it follows from Theorem 3.6 there there exists another $(q+1)$-tuple
$(N', m_1', \cdots, m_q') \in\N^{q+1}$, such that the inequalities (\ref{4.27})-(\ref{4.31'}) hold for any
$1\le k\le q$. By the similar arguments in Step A, we obtain the corresponding numbers $N_{+}^{'e}$ defined in
the same way as in (\ref{4.33}), and it yields
$$  N_+^{'e} \ge \frac{n}{2}. $$

Then the symmetry of $\chi$ and $\hat{\chi}$ yields $N_-^e = N_+^{'e}$. Thus we obtain
$$  q \ge N_+^e + N_-^e\ge \frac{n}{2} + \frac{n}{2} = n.  $$
That is, the total number of distinct closed characteristics on such $\Sigma$ in Theorem 1.2 is at least
$n$.

{\bf Step C.} Then the non-hyperbolicity of these closed characteristics can be easily obtained according to their
even index information given in the definition (\ref{4.21})-(\ref{4.23}) of $N_+^e$ and $N_-^e$, because if $y_k$
is hyperbolic, it must satisfy $i(y_k^{2m_k})$ being odd in the case of $n$ being even.

\medskip

Based on Theorem 3.6 below, we can generalize Theorem 1.2 about closed characteristics on star-shaped hypersurfaces to
closed Reeb orbits of contact forms on a broad class of prequantization bundles, which improves Theorem 2.1 and
Theorem 2.10 of \cite{GGM}. To clarify it, we review some terminologies from contact geometry, following Section 2
of \cite{GGM}.

Let $(M^{2n+1}, \xi)$ be a closed contact manifold satisfying $c_1(\xi)|_{\pi_2(M)} = 0$ and $\alpha$ be the contact form
supporting the contact structure $\xi$. A non-degenerate periodic orbit $\gamma$ is called {\it good} if its Conley-Zehnder
index $\mu(\gamma)$ has the same parity as that of the underlying simple closed orbit. Note that the Maslov-type index of
$\gamma$ equals to $\mu(\gamma)-1$. In the following, let
$(M^{2n+1}, \xi)$ be a prequantization circle bundle over closed integral symplectic manifolds
$(B^{2n}, \omega)$, i.e., the first Chern class of the principle bundle $M \rightarrow B$ is $-[\omega]$.
Denote by $\chi(B)$ the Euler characteristic of $B$ and by
\[c_B :=\mathrm{ inf}\{k \in\N ~| ~\exists S \in  \pi_2(B) ~\mathrm{with}~ \langle c_1(TB), S\rangle     = k\}\]
its minimal Chern number. We impose the following condition which is weaker than the condition (F) in Section 2 of \cite{GGM}:
\medskip

($\mathbf{A}$) (i) The manifold $(M^{2n+1}, \xi)$ admits a strong symplectic filling $(W, \Omega)$ which is symplectically
aspherical, i.e., $\Omega|_{\pi_2(W)} = 0$ and $c_1(TW)|_{\pi_2(W)} = 0$, and the map $\pi_1(M)\rightarrow \pi_1(W)$ induced
by the inclusion is injective.

(ii) The contact form $\alpha$ is non-degenerate, the mean index $\hat{\mu}(\gamma)$ is nonzero for every contractible
periodic orbit $\gamma$ of $\alpha$ and has no contractible good periodic orbits $\gamma$ such that $\mu(\gamma)=0$ if $n$
is odd or $\mu(\gamma)\in \{0,\pm 1\}$ if $n$ is even.

\medskip

{\bf Theorem 1.3.} {\it Let $(M^{2n+1}, \xi)$ be a prequantization $S^1$-bundle of a closed symplectic
manifold $(B, \omega)$ such that $\omega|_{\pi_2(B)}= 0$ and $c_B > n/2$ and, furthermore,
$H_k(B;\Q) = 0$ for every odd $k$ or $c_B > n$. Let $\alpha$ be a contact form supporting
$\xi$ and assume that $M$ and $\alpha$ satisfy condition $(\mathbf{A})$. Then $\alpha$ carries at
least $r_B$ geometrically distinct contractible periodic orbits. Furthermore, if there exist finitely many geometrically
distinct contractible closed orbits, then $\alpha$ carries at least $r^{non-hyp}_B$ geometrically distinct contractible
non-hyperbolic periodic orbits, where $r^{non-hyp}_B:=r_B-\mathrm{dim}H_n(B;\Q)$ and
\be r_B:=\left\{\begin{array}{ll}
    \chi(B)+ 2~\mathrm{dim}H_n(B;\Q), &\quad {\it if}\;\;n ~\mathrm{is}~ \mathrm{odd}, \\
     \chi(B)+ 4~\mathrm{dim}H_{n-1}(B;\Q), &\quad {\it if}\;\;n ~\mathrm{is}~ \mathrm{even}.\end{array}\right. \nn\ee}

\medskip

{\bf Remark 1.4.} (1) The proof of Theorem 1.3 follows the proofs of Theorems 2.1 and 2.10 of \cite{GGM} via replacing the enhanced
common jump theorem of H. Duan, Y. Long and W. Wang by our Theorem 3.6 and the proof of Theorem 1.2 below. We should emphasize that
under nonzero mean index condition, there are similar Morse inequalities in the setting of equivariant symplectic homology to
(\ref{2.16})-(\ref{2.17}) below, cf., (68) in p.221 of \cite{HM}.

(2) Note that we are unable to weaken the condition (NF) in Section 2 of \cite{GGM} as ($\mathbf{A}$), since the proof in \cite{GGM}
relies on the machinery of positive equivariant symplectic homology. A remarkable observation by F. Bourgeois and A. Oancea in
Section 4.1.2 of \cite{BO} is that under suitable additional assumptions that all closed contractible Reeb orbits on $M$ are
non-degenerate and have Conley-Zehnder index strictly greater than $3-n$, the positive equivariant symplectic homology is defined
and well-defined even when $M$ does not have a symplectic filling. But the existence of closed contractible Reeb orbit with negative
mean index will destroy this assumption and then the positive equivariant symplectic homology is not well defined, thus we are unable to
weaken the index-positive condition in (NF) of \cite{GGM} to the nonzero mean index condition.

\medskip

An important result of V. Bangert and W. Klingenberg in \cite{BK1983} implies that if $c$ is a closed geodesic on a compact Riemannian
(or Finsler) manifold $M$ such that it possesses zero mean index and $c$ is neither homologically invisible nor an absolute minimum of
the energy functional, then there exist infinitely many closed geodesics on $M$. In fact, we tend to believe that when the number of
prime closed geodesics on a compact Finsler manifold or prime closed characteristics on a compact star-shaped hypersurface in
$\R^{2n}$ is finite, then each one of them must be homologically visible as well as variationally visible (cf. \cite{BK1983} and
\cite{Lon4} for definitions). Motivated by the result in \cite{BK1983} and the well-known weakly non-resonant ellipsoids, we tend to
believe that the following conjecture for closed characteristics should hold:

\medskip

{\bf Conjecture 1.5.} {\it If there exist only finitely many geometrically distinct closed characteristics on a compact star-shaped
hypersurface $\Sg$ in $\R^{2n}$, then no one of them possesses zero mean index.}

\medskip

Our third goal of this paper is to give a positive answer to this conjecture below in the case of $n=3$ for
non-degenerate star-shaped hypersurfaces. But up to our knowledge, this conjecture seems to be rather challenging in its
full generality.

\medskip

{\bf Theorem 1.6.} {\it If there exist only finitely many geometrically distinct closed characteristics on a compact non-degenerate
star-shaped hypersurface $\Sg$ in $\R^{6}$, then every prime closed characteristic possesses nonzero mean
index.}
\medskip

Using Theorem 1.6, we can remove the nonzero mean index condition in Theorem 1.2 in the case of $n=3$:

\medskip

{\bf Corollary 1.7.} {\it If $\Sigma$ is a compact non-degenerate perfect star-shaped hypersurface in $\R^{6}$,
then there exist at least three geometrically distinct closed characteristics. Furthermore, if the total number of prime
closed characteristics on $\Sigma$ is finite, then $\Sigma$ carries at least three geometrically distinct closed
characteristics with odd Maslov-type indices and at least two of them are non-hyperbolic.}

\medskip

Note that one may generalize our Theorem 1.6 to contact manifolds by the idea of our proof of Theorem 1.6, especially the
key observation in our Lemma 5.1 that the Viterbo index $i(y^m)$ always equals to $-4$ for every closed characteristic
$(\tau, y)$ with zero mean index and all $m\in\N$.

\medskip

In this paper, let $\N$, $\N_0$, $\Z$, $\Q$, $\R$, $\C$ and $\R^+$ denote the sets of natural integers,
non-negative integers, integers, rational numbers, real numbers, complex numbers and positive real
numbers respectively. We define the functions
\bea [a]&=&\max{\{k\in {\bf Z}\mid k\leq a\}},\qquad  \{a\}=a-[a],\\
E(a)&=&\min{\{k\in{\bf Z}\mid k\geq a\}},\  \varphi(a)=E(a)-[a].\lb{1.2}\eea
Denote by $a\cdot b$ and $|a|$ the standard inner product and norm in $\R^{2n}$. Denote by $\langle\cdot,\cdot\rangle$ and $\|\cdot\|$
the standard $L^2$ inner product and $L^2$ norm. For an $S^1$-space $X$, we denote by
$X_{S^1}$ the homotopy quotient of $X$ by $S^1$, i.e., $X_{S^1}=S^\infty\times_{S^1}X$,
where $S^\infty$ is the unit sphere in an infinite dimensional {\it complex} Hilbert space.
In this paper we use $\Q$ coefficients for all homological and cohomological modules. By $t\to a^+$, we
mean $t>a$ and $t\to a$.

This paper is organized as follows. In section 2, we recall some variational settings, which includes the $S^1$-critical module, mean index identities for closed characteristics and the Morse inequality. In section 3, based on \cite{LoZ} and \cite{DLW}, we establish the {\bf generalized common index jump theorem} (see Theorem 3.6 below) for symplectic paths. Then the proofs of Theorem 1.2, Theorem 1.3 and Theorem 1.6 are given in Section 4 and Section 5 respectively.

\setcounter{figure}{0}
\setcounter{equation}{0}
\section{Mean index identities for closed characteristics on compact star-shaped hypersurfaces in $\R^{2n}$}

In this section, we briefly review the mean index identities for
closed characteristics on $\Sg\in\H_{st}(2n)$ developed in \cite{LLW} which will be needed in Section 4. All
the details of proofs can be found in \cite{LLW}.
Now we fix a $\Sg\in\H_{st}(2n)$ and assume the following condition on $\T(\Sg)$:

\medskip

(F) {\bf There exist only finitely many geometrically distinct prime closed characteristics
$\{(\tau_j, y_j)\}_{1\le j\le k}$ on $\Sigma$.}

\medskip

Let $j: \R^{2n}\rightarrow\R$ be the gauge function of $\Sigma$, i.e., $j(\lambda x)=\lambda$ for $x\in\Sigma$
and $\lambda\ge0$, then $j\in C^3(\R^{2n}\bs\{0\}, \R)\cap C^0(\R^{2n}, \R)$ and $\Sigma=j^{-1}(1)$.
Let $\hat{\tau}=\inf_{1\leq j\leq k}{\tau_j}$ and $T$ be a fixed positive constant. Then following \cite{Vit1}
and Section 2 of \cite{LLW}, for any $a>\frac{\hat{\tau}}{T}$, we can construct a function
$\varphi_a\in C^{\infty}({\bf R}, {\bf R}^+)$ which has $0$ as its unique critical point in $[0, +\infty)$. Moreover,
$\frac{\varphi_a^{\prime}(t)}{t}$ is strictly decreasing for $t>0$ together with $\varphi_a(0)=0=\varphi_a^{\prime}(0)$
and $\varphi_a^{\prime\prime}(0)=1=\lim_{t\rightarrow 0^+}\frac{\varphi_a^{\prime}(t)}{t}$. The precise definition of
$\varphi_a$ and the dependence of $\varphi_a$ on $a$ are given in Lemma 2.2 and Remark 2.3 of \cite{LLW} respectively.
As in \cite{LLW}, we define a Hamiltonian function
$H_a\in C^{3}({\bf R}^{2n} \setminus\{0\},{\bf R})\cap C^{1}({\bf R}^{2n},{\bf R})$ satisfying $H_a(x)=a\vf_a(j(x))$
on $U_A=\{x\mid a\vf_a(j(x))\leq A\}$ for some large $A$, and $H_a(x)=\frac{1}{2}\ep_a|x|^2$ outside some even larger
ball with $\ep_a>0$ small enough such that outside $U_A$ both $\nabla H_a(x)\not= 0$ and $H_a^{\prime\prime}(x)<\ep_a$
hold as in Lemmas 2.2, 2.4 and Proposition 2.5 of \cite{LLW}.

We consider the following fixed period problem
\be  \dot{x}(t)=JH_a^\prime(x(t)),\quad x(0)=x(T).  \lb{2.1}\ee
Then solutions of (\ref{2.1}) are $x\equiv 0$ and $x=\rho y(\tau t/T)$ with
$\frac{\vf_a^\prime(\rho)}{\rho}=\frac{\tau}{aT}$, where $(\tau, y)$ is a solution of (\ref{1.1}). In particular,
non-zero solutions of (\ref{2.1}) are in one to one correspondence with solutions of (\ref{1.1}) with period
$\tau<aT$.

For any $a>\frac{\hat{\tau}}{T}$, we can choose some large constant $K=K(a)$ such that
\be H_{a,K}(x) = H_a(x)+\frac{1}{2}K|x|^2   \lb{2.2}\ee
is a strictly convex function, that is,
\be (\nabla H_{a, K}(x)-\nabla H_{a, K}(y), x-y) \geq \frac{\ep}{2}|x-y|^2,  \lb{2.3}\ee
for all $x, y\in {\bf R}^{2n}$, and some positive $\ep$. Let $H_{a,K}^*$ be the Fenchel dual of $H_{a,K}$
defined by
$$  H_{a,K}^\ast (y) = \sup\{x\cdot y-H_{a,K}(x)\;|\; x\in \R^{2n}\}.   $$
The dual action functional on $X=W^{1, 2}({\bf R}/{T {\bf Z}}, {\bf R}^{2n})$ is defined by
\be F_{a,K}(x) = \int_0^T{\left[\frac{1}{2}(J\dot{x}-K x,x)+H_{a,K}^*(-J\dot{x}+K x)\right]dt}. \lb{2.4}\ee
Then $F_{a,K}\in C^{1,1}(X, \R)$ and for $KT\not\in 2\pi{\bf Z}$, $F_{a,K}$ satisfies the
Palais-Smale condition and $x$ is a critical point of $F_{a, K}$ if and only if it is a solution of (\ref{2.1}). Moreover,
$F_{a, K}(x_a)<0$ and it is independent of $K$ for every critical point $x_a\neq 0$ of $F_{a, K}$.

When $KT\notin 2\pi{\bf Z}$, the map $x\mapsto -J\dot{x}+Kx$ is a Hilbert space isomorphism between
$X=W^{1, 2}({\bf R}/({T {\bf Z}}); {\bf R}^{2n})$ and $E=L^{2}({\bf R}/(T {\bf Z}),{\bf R}^{2n})$. We denote its inverse
by $M_K$ and the functional
\be \Psi_{a,K}(u)=\int_0^T{\left[-\frac{1}{2}(M_{K}u, u)+H_{a,K}^*(u)\right]dt}, \qquad \forall\,u\in E. \lb{2.5}\ee
Then $x\in X$ is a critical point of $F_{a,K}$ if and only if $u=-J\dot{x}+Kx$ is a critical point of $\Psi_{a, K}$.

Suppose $u$ is a nonzero critical point of $\Psi_{a, K}$.
Then the formal Hessian of $\Psi_{a, K}$ at $u$ is defined by
\be Q_{a,K}(v)=\int_0^T(-M_K v\cdot v+H_{a,K}^{*\prime\prime}(u)v\cdot v)dt,  \lb{2.6}\ee
which defines an orthogonal splitting $E=E_-\oplus E_0\oplus E_+$ of $E$ into negative, zero and positive subspaces.
The index and nullity of $u$ are defined by $i_K(u)=\dim E_-$ and $\nu_K(u)=\dim E_0$ respectively.
Similarly, we define the index and nullity of $x=M_Ku$ for $F_{a, K}$, we denote them by $i_K(x)$ and
$\nu_K(x)$. Then we have
\be  i_K(u)=i_K(x),\quad \nu_K(u)=\nu_K(x),  \lb{2.7}\ee
which follow from the definitions (\ref{2.4}) and (\ref{2.5}). The following important formula was proved in
Lemma 6.4 of \cite{Vit2}:
\be  i_K(x) = 2n([KT/{2\pi}]+1)+i^v(x) \equiv d(K)+i^v(x),   \lb{2.8}\ee
where the Viterbo index $i^v(x)$ does not depend on K, but only on $H_a$.

By the proof of Proposition 2 of \cite{Vit1}, we have that $v\in E$ belongs to the null space of $Q_{a, K}$
if and only if $z=M_K v$ is a solution of the linearized system
\be  \dot{z}(t) = JH_a''(x(t))z(t).  \lb{2.9}\ee
Thus the nullity in (\ref{2.7}) is independent of $K$, which we denote by $\nu^v(x)\equiv \nu_K(u)= \nu_K(x)$.

By Proposition 2.11 of \cite{LLW}, the index $i^v(x)$ and nullity $\nu^v(x)$ coincide with those defined for
the Hamiltonian $H(x)=j(x)^\alpha$ for all $x\in\R^{2n}$ and some $\aa\in (1,2)$. Especially
$1\le \nu^v(x)\le 2n-1$ always holds.

For every closed characteristic $(\tau, y)$ on $\Sigma$, let $aT>\tau$ and choose $\vf_a$ as above.
Determine $\rho$ uniquely by $\frac{\vf_a'(\rho)}{\rho}=\frac{\tau}{aT}$. Let $x=\rho y(\frac{\tau t}{T})$.
Then we define the index $i(\tau,y)$ and nullity $\nu(\tau,y)$ of $(\tau,y)$ by
$$ i(\tau,y)=i^v(x), \qquad \nu(\tau,y)=\nu^v(x). $$
Then the mean index of $(\tau,y)$ is defined by
\bea \hat i(\tau,y) = \lim_{m\rightarrow\infty}\frac{i(m\tau,y)}{m}.  \nn\eea
Note that by Proposition 2.11 of \cite{LLW}, the index and nullity are well defined and are independent of the
choice of $a$. For a closed characteristic $(\tau,y)$ on $\Sigma$, we simply denote by $y^m\equiv(m\tau,y)$
the m-th iteration of $y$ for $m\in\N$.

We have a natural $S^1$-action on $X$ or $E$ defined by
$$  \theta\cdot u(t)=u(\theta+t),\quad\forall\, \theta\in S^1, \, t\in\R.  $$
Clearly both of $F_{a, K}$ and $\Psi_{a, K}$ are $S^1$-invariant. For any $\kappa\in\R$, we denote by
\bea
\Lambda_{a, K}^\kappa &=& \{u\in L^{2}({\bf R}/({T {\bf Z}}); {\bf R}^{2n})\;|\;\Psi_{a,K}(u)\le\kappa\},  \nn\\
X_{a, K}^\kappa &=& \{x\in W^{1, 2}({\bf R}/(T {\bf Z}),{\bf R}^{2n})\;|\;F_{a, K}(x)\le\kappa\}.  \nn\eea
For a critical point $u$ of $\Psi_{a, K}$ and the corresponding $x=M_K u$ of $F_{a, K}$, let
\bea
\Lm_{a,K}(u) &=& \Lm_{a,K}^{\Psi_{a, K}(u)}
   = \{w\in L^{2}(\R/(T\Z), \R^{2n}) \;|\; \Psi_{a, K}(w)\le\Psi_{a,K}(u)\},  \nn\\
X_{a,K}(x) &=& X_{a,K}^{F_{a,K}(x)} = \{y\in W^{1, 2}(\R/(T\Z), \R^{2n}) \;|\; F_{a,K}(y)\le F_{a,K}(x)\}. \nn\eea
Clearly, both sets are $S^1$-invariant. Denote by $\crit(\Psi_{a, K})$ the set of critical points of $\Psi_{a, K}$.
Because $\Psi_{a,K}$ is $S^1$-invariant, $S^1\cdot u$ becomes a critical orbit if $u\in \crit(\Psi_{a, K})$.
Note that by the condition (F), the number of critical orbits of $\Psi_{a, K}$
is finite. Hence as usual we can make the following definition.

\medskip

{\bf Definition 2.1.} {\it Suppose $u$ is a nonzero critical point of $\Psi_{a, K}$, and $\Nn$ is an $S^1$-invariant
open neighborhood of $S^1\cdot u$ such that $\crit(\Psi_{a,K})\cap (\Lm_{a,K}(u)\cap \Nn) = S^1\cdot u$.
Then the $S^1$-critical module of $S^1\cdot u$ is defined by
$$ C_{S^1,\; q}(\Psi_{a, K}, \;S^1\cdot u)
=H_{q}((\Lambda_{a, K}(u)\cap\Nn)_{S^1},\; ((\Lambda_{a,K}(u)\setminus S^1\cdot u)\cap\Nn)_{S^1}). $$
Similarly, we define the $S^1$-critical module $C_{S^1,\; q}(F_{a, K}, \;S^1\cdot x)$ of $S^1\cdot x$
for $F_{a, K}$.}

\medskip

We fix $a$ and let $u_K\neq 0$ be a critical point of $\Psi_{a, K}$ with multiplicity $\mul(u_K)=m$,
that is, $u_K$ corresponds to a closed characteristic $(\tau, y)\subset\Sigma$ with $(\tau, y)$
being $m$-iteration of
some prime closed characteristic. Precisely, we have $u_K=-J\dot x+Kx$ with $x$
being a solution of (\ref{2.1}) and $x=\rho y(\frac{\tau t}{T})$ with
$\frac{\vf_a^\prime(\rho)}{\rho}=\frac{\tau}{aT}$.
Moreover, $(\tau, y)$ is a closed characteristic on $\Sigma$ with minimal period $\frac{\tau}{m}$.
For any $p\in\N$ satisfying $p\tau<aT$, we choose $K$
such that $pK\notin \frac{2\pi}{T}\Z$, then the $p$th iteration $u_{pK}^p$ of $u_K$ is given by $-J\dot x^p+pKx^p$,
where $x^p$ is the unique solution of (\ref{2.1}) corresponding to $(p\tau, y)$
and is a critical point of $F_{a, pK}$, that
is, $u_{pK}^p$ is the critical point of $\Psi_{a, pK}$ corresponding to $x^p$.

\medskip

{\bf Lemma 2.2.} (cf. Proposition 4.2 and Remark 4.4 of \cite{LLW} ) {\it If $u_{pK}^p$ is non-degenerate,
i.e., $\nu_{pK}(u_{pK}^p)=1$, let
$\bb(x^p)=(-1)^{i_{pK}(u_{pK}^p)-i_{K}(u_{K})}=(-1)^{i^v(x^p)-i^v(x)}$, then}
\bea &&C_{S^1,q-d(pK)+d(K)}(F_{a,K},S^1\cdot x^p)\nn\\
&&\qquad= C_{S^1,q}(F_{a,pK},S^1\cdot x^p)\nn\\
&&\qquad= C_{S^1,q}(\Psi_{a,pK},S^1\cdot u^p_{pK}) \nn\\
&&\qquad= \left\{\begin{array}{ll}
     \Q, &\quad {\it if}\;\; q=i_{pK}(u_{pK}^p)\;\;{\it and}\;\;\bb(x^p)=1, \\
     0, &\quad {\it otherwise}. \end{array}\right.  \lb{2.10}\eea

\medskip

{\bf Theorem 2.3.} (cf. Theorem 1.1 of \cite{LLW} and Theorem 1.2 of \cite{Vit2}) {\it Suppose that
$\Sg\in\H_{st}(2n)$ satisfying $^\#\T(\Sg)<+\infty$. Denote by $\{(\tau_j,y_j)\}_{1\le j\le k}$ all the
geometrically distinct prime closed characteristics. Then the following identities hold
\bea \sum_{1\le j\le k \atop \hat{i}(y_j)>0}\frac{\hat{\chi}(y_j)}{\hat{i}(y_j)}=\frac{1}{2},\qquad
\sum_{1\le j\le k \atop \hat{i}(y_j)<0}\frac{\hat{\chi}(y_j)}{\hat{i}(y_j)}=0,\lb{2.11}\eea
where $\hat{\chi}(y)\in\Q$ is the average Euler characteristic given by Definition 4.8 and Remark 4.9 of \cite{LLW}.

In particular, if all $y^m$'s are non-degenerate for $m\ge 1$, then
\bea \hat{\chi}(y)=\left\{\begin{array}{ll}
     (-1)^{i(y)}, &\quad {\it if}\;\; i(y^2)-i(y)\in 2\Z, \\
     \frac{(-1)^{i(y)}}{2}, &\quad {\it otherwise}. \end{array}\right.  \lb{2.12}\eea}

\medskip

Let $F_{a, K}$ be a functional defined by (\ref{2.4}) for some $a, K\in\R$ large enough and let $\ep>0$ be
small enough such that $[-\ep, 0)$ contains no critical values of $F_{a, K}$. For $b$ large enough,
The normalized Morse series of $F_{a, K}$ in $ X^{-\ep}\setminus X^{-b}$
is defined, as usual, by
\be  M_a(t)=\sum_{q\ge 0,\;1\le j\le p} \dim C_{S^1,\;q}(F_{a, K}, \;S^1\cdot v_j)t^{q-d(K)},  \lb{2.13}\ee
where we denote by $\{S^1\cdot v_1, \ldots, S^1\cdot v_p\}$ the critical orbits of $F_{a, K}$ with critical
values less than $-\ep$. The Poincar\'e series of $H_{S^1, *}( X, X^{-\ep})$ is $t^{d(K)}Q_a(t)$, according
to Theorem 5.1 of \cite{LLW}, if we set $Q_a(t)=\sum_{k\in \Z}{q_kt^k}$, then
$$   q_k=0, \quad \forall\;k\in \mathring {I},  $$
where $I$ is an interval of $\Z$ such that $I \cap [i(\tau, y), i(\tau, y)+\nu(\tau, y)-1]=\emptyset$ for all
closed characteristics $(\tau,\, y)$ on $\Sigma$ with $\tau\ge aT$. Then by Section 6 of \cite{LLW}, we have
$$  M_a(t)-\frac{1}{1-t^2}+Q_a(t) = (1+t)U_a(t),   $$
where $U_a(t)=\sum_{i\in \Z}{u_it^i}$ is a Laurent series with nonnegative coefficients.
If there is no closed characteristic with $\hat{i}=0$, then
\be   M(t)-\frac{1}{1-t^2}=(1+t)U(t),    \lb{2.14}\ee
where $M(t)=\sum_{p\in \Z}{M_pt^p}$ denotes $M_a(t)$ as $a$ tends to infinity. In addition, we also denote by
$b_p$ the coefficient of $t^p$ of $\frac{1}{1-t^2}=\sum_{p\in \Z}{b_pt^p}$, i.e. there holds $b_p=1$, for all
$p\in 2\N_0$, and $b_p=0$ for all $p \not\in 2\N_0$.

For any two positive integers $n_1$ and $n_2$, it follows from (\ref{2.14}) that
\bea &&\sum_{p=-2n_1+1}^{2n_2+1}{M_pt^p}-\sum_{p=-2n_1+1}^{2n_2+1}b_pt^p\nn\\
           &&\qquad=  (1+t)\sum_{p=-2n_1}^{2n_2+1}u_pt^p
           -u_{2n_2+1}t^{2n_2+2}-u_{-2n_1}t^{-2n_1}, \eea
which, through letting $t=-1$, yields the following Morse inequality
\bea \sum_{p=-2n_1+1}^{2n_2+1}{(-1)^pM_p} \le \sum_{p=-2n_1+1}^{2n_2+1}(-1)^pb_p.\lb{2.16}\eea

Similarly we have
\bea \sum_{p=-2n_1}^{2n_2}{(-1)^pM_p} \ge \sum_{p=-2n_1}^{2n_2}(-1)^pb_p.\lb{2.17}\eea

\setcounter{figure}{0}
\setcounter{equation}{0}
\section{The generalized common index jump theorem for symplectic paths}

In \cite{Lon2} of 1999, Y. Long established the basic normal form
decomposition of symplectic matrices. Based on this result he
further established the precise iteration formulae of indices of
symplectic paths in \cite{Lon3} of 2000.

As in \cite{Lon3}, denote by
\bea
N_1(\lm, b) &=& \left(\begin{array}{ll}\lm & b\\
                                0 & \lm \end{array}\right), \qquad {\rm for\;}\lm=\pm 1, \; b\in\R, \lb{3.1}\\
D(\lm) &=& \left(\begin{array}{ll}\lm & 0\\
                      0 & \lm^{-1} \end{array}\right), \qquad {\rm for\;}\lm\in\R\bs\{0, \pm 1\}, \lb{3.2}\\
R(\th) &=& \left(\begin{array}{ll}\cos\th & -\sin\th \\
                           \sin\th & \cos\th \end{array}\right), \qquad {\rm for\;}\th\in (0,\pi)\cup (\pi,2\pi), \lb{3.3}\\
N_2(e^{\th\sqrt{-1}}, B) &=& \left(\begin{array}{ll} R(\th) & B \\
                  0 & R(\th) \end{array} \right), \qquad {\rm for\;}\th\in (0,\pi)\cup (\pi,2\pi)\;\; {\rm and}\; \nn\\
        && \quad B=\left(\begin{array}{lll} b_1 & b_2\\
                                  b_3 & b_4 \end{array}\right)\; {\rm with}\; b_j\in\R, \;\;
                                         {\rm and}\;\; b_2\not= b_3. \lb{3.4}\eea
Here $N_2(e^{\th\sqrt{-1}}, B)$ is non-trivial if $(b_2-b_3)\sin\theta<0$, and trivial
if $(b_2-b_3)\sin\theta>0$.

As in \cite{Lon3}, the $\diamond$-sum (direct sum) of any two real matrices is defined by
$$ \left(\begin{array}{lll}A_1 & B_1\\ C_1 & D_1 \end{array}\right)_{2i\times 2i}\diamond
      \left(\begin{array}{lll}A_2 & B_2\\ C_2 & D_2 \end{array}\right)_{2j\times 2j}
=\left(\begin{array}{llll}A_1 & 0 & B_1 & 0 \\
                                   0 & A_2 & 0& B_2\\
                                   C_1 & 0 & D_1 & 0 \\
                                   0 & C_2 & 0 & D_2\end{array}\right). $$

For every $M\in\Sp(2n)$, the homotopy set $\Omega(M)$ of $M$ in $\Sp(2n)$ is defined by
$$ \Om(M)=\{N\in\Sp(2n)\,|\,\sg(N)\cap\U=\sg(M)\cap\U\equiv\Gamma,
                    \;\nu_{\om}(N)=\nu_{\om}(M),\, \forall\om\in\Gamma\}, $$
where $\sg(M)$ denotes the spectrum of $M$,
$\nu_{\om}(M)\equiv\dim_{\C}\ker_{\C}(M-\om I)$ for $\om\in\U$.
The component $\Om^0(M)$ of $P$ in $\Sp(2n)$ is defined by
the path connected component of $\Om(M)$ containing $M$.

\medskip

{\bf Lemma 3.1.} (cf. \cite{Lon3}, Lemma 9.1.5 and List 9.1.12 of \cite{Lon4})
{\it For $M\in\Sp(2n)$ and $\om\in\U$, the splitting number $S_M^\pm(\om)$ (cf. Definition 9.1.4 of \cite{Lon4}) satisfies
\begin{eqnarray}
S_M^{\pm}(\om) &=& 0, \qquad {\rm if}\;\;\om\not\in\sg(M).  \lb{3.5}\\
S_{N_1(1,a)}^+(1) &=& \left\{\begin{array}{lll}1, &\quad {\rm if}\;\; a\ge 0, \\
0, &\quad {\rm if}\;\; a< 0. \end{array}\right. \lb{3.6}\eea

For any $M_i\in\Sp(2n_i)$ with $i=0$ and $1$, there holds }
\be S^{\pm}_{M_0\diamond M_1}(\om) = S^{\pm}_{M_0}(\om) + S^{\pm}_{M_1}(\om),
    \qquad \forall\;\om\in\U. \lb{3.7}\ee

We have the following decomposition theorem

\medskip

{\bf Theorem 3.2.} (cf. \cite{Lon3} and Theorem 1.8.10 of \cite{Lon4}) {\it For
any $M\in\Sp(2n)$, there is a path $f:[0,1]\to\Om^0(M)$ such that $f(0)=M$ and
\be f(1) = M_1\diamond\cdots\diamond M_k,  \lb{3.8}\ee
where each $M_i$ is a basic normal form listed in (\ref{3.1})-(\ref{3.4})
for $1\leq i\leq k$.}

\medskip

For every $\ga\in\mathcal{P}_\tau(2n)\equiv\{\ga\in C([0,\tau],Sp(2n))\ |\ \ga(0)=I_{2n}\}$, we extend
$\ga(t)$ to $t\in [0,m\tau]$ for every $m\in\N$ by
\bea \ga^m(t)=\ga(t-j\tau)\ga(\tau)^j \quad \forall\;j\tau\le t\le (j+1)\tau \;\;
               {\rm and}\;\;j=0, 1, \ldots, m-1, \lb{3.9}\eea
as in p.114 of \cite{Lon2}. As in \cite{LoZ} and \cite{Lon4}, we denote the Maslov-type indices of
$\ga^m$ by $(i(\ga,m),\nu(\ga,m))$.

Then the following iteration formula from \cite{LoZ} and \cite{Lon4} can be obtained.

\medskip

{\bf Theorem 3.3.} (cf. Theorem 9.3.1 of \cite{Lon4}) {\it For any path $\ga\in\mathcal{P}_\tau(2n)$,
let $M=\ga(\tau)$ and $C(M)=\sum_{0<\th<2\pi}S_M^-(e^{\sqrt{-1}\th})$. We extend $\ga$ to $[0,+\infty)$
by its iterates. Then for any $m\in\N$ we have
\bea &&i(\ga,m)
= m(i(\ga,1)+S^+_{M}(1)-C(M))\nn\\
&&\qquad+ 2\sum_{\th\in(0,2\pi)}E\left(\frac{m\th}{2\pi}\right)S^-_{M}(e^{\sqrt{-1}\th}) - (S_M^+(1)+C(M)) \lb{3.10}\eea
and
\be \hat{i}(\ga,1) = i(\ga,1) + S^+_{M}(1) - C(M) + \sum_{\th\in(0,2\pi)}\frac{\th}{\pi}S^-_{M}(e^{\sqrt{-1}\th}). \lb{3.11}\ee}

\medskip

{\bf Theorem 3.4.} {\it Fix an integer $q>0$. Let $\mu_i\ge 0$ and $\bb_i$ be integers for all $i=1,\cdots,q$. Let $\aa_{i,j}$
be positive numbers for $j=1,\cdots,\mu_i$ and $i=1,\cdots,q$. Let $\dl\in(0,\frac{1}{2})$ satisfying
$\dl\max\limits_{1\le i\le q}\mu_i<\frac{1}{2}$. Suppose $D_i \equiv \bb_i+\sum\limits_{j=1}^{\mu_i}\aa_{i,j}\neq 0$ for
$i=1,\cdots,q$. Then there exist infinitely many $(N, m_1,\cdots,m_q)\in\N^{q+1}$ such that
\bea
&& m_i\bb_i+\sum_{j=1}^{\mu_i}E(m_i\aa_{i,j}) =
      \varrho_i N+\Delta_i, \qquad \forall\ 1\le i\le q.  \lb{3.12}\\
&& \min\{\{m_i\aa_{i,j}\}, 1-\{m_i\aa_{i,j}\}\} < \dl,\  \forall\ j=1,\cdots,\mu_i, 1\le i\le q, \lb{3.13}\\
&& m_i\aa_{i,j}\in\N,\  {\rm if} \  \aa_{i,j}\in\Q,   \lb{3.14}\eea
where
\bea \varrho_i=\left\{\begin{array}{cc}1, &{\rm if}\ D_i>0, \cr
                                     -1, &{\rm if}\  D_i<0, \end{array}\right.\quad \Delta_i=\sum_{0<\{m_i\aa_{i,j}\}<\dl}1,\quad \forall\ 1\le i\le q.\lb{3.15}\eea}

\medskip

{\bf Remark 3.5.} When $D_i>0$ for all $1\le i\le q$, this is precisely the Theorem 4.1 of \cite{LoZ} (also cf.
Theorem 11.1.1 of \cite{Lon4}).

\medskip

{\bf Proof of Theorem 3.4.}  By assumption $D_i\neq 0$, $\forall\ 1\le i\le q$, we further assume that there exists some
integer $0\le q_0\le q$ with $D_i<0$ for $0\le i\le q_0$ and $D_i>0$ for $q_0+1\le i\le q$. Next we will do with both of
these two cases simultaneously. In fact we only need to use $\varrho_i N$ and
$m_i=\left(\left[\frac{\varrho_i N}{M D_i}\right]+\chi_i\right)M$ to replace the corresponding $N$ and $m_i$s in the
proof of Theorem 4.1 of \cite{LoZ} (cf. Theorem 11.1.1 of \cite{Lon4}). For reader's convenience and because the proof is
almost self-contained, in the following we only give some different points and details.

In order to get (\ref{3.12}), we consider
\bea m_i D_i&=&\frac{\varrho_i N}{MD_i}MD_i-\left\{\frac{\varrho_i N}{M D_i}\right\}MD_i+\chi_i MD_i\nn\\
&=&\varrho_i N+\left(\chi_i-\left\{\frac{\varrho_i N}{M D_i}\right\}\right)MD_i,\quad \forall\ 1\le i\le q,\lb{3.16}\eea
where, to get (\ref{3.14}), we require $M\in\N$ to satisfy $M\aa_{i,j}\in\N$ when $\aa_{i,j}\in\Q$ for $j=1,\cdots,\mu_i$,
and $\chi_i\in\{0,1\}$ will be determined later.

Set \bea m_i=\left(\left[\frac{\varrho_i N}{M D_i}\right]+\chi_i\right)M.\lb{3.17}\eea

Then by (\ref{1.2}) and (\ref{3.16}), following the proofs from (4.11) and (4.13) of \cite{LoZ} (or (11.1.11) to (11.1.13) of
\cite{Lon4}) and using $\Delta_i$ and $\delta$ defined there, it yields
\bea &&m_i\bb_i+\sum_{j=1}^{\mu_i}E(m_i\aa_{i,j}) \nn\\
&&\quad= m_i D_i+\sum_{j=1}^{\mu_i}(\varphi(m_i\alpha_{i,j})-\{m_i\alpha_{i,j}\}) \nn\\
&&\quad= \varrho_i N+\left(\chi_i-\left\{\frac{\varrho_i N}{M D_i}\right\}\right)MD_i\nn\\
&&\qquad\ +\sum_{j=1}^{\mu_i}(\varphi(m_i\alpha_{i,j})-\{m_i\alpha_{i,j}\})\nn\\
&&\quad= \varrho_i N+\left(\chi_i-\left\{\frac{\varrho_i N}{M D_i}\right\}\right)MD_i+\Delta_i \nn\\
  &&\qquad\ -\sum_{0<\{m_i\alpha_{i,j}\}<\dl}\{m_i\alpha_{i,j}\}+\sum_{0<1-\{m_i\aa_{i,j}\}<\delta}(1-\{m_i\alpha_{i,j}\}),\lb{3.18}\eea
which, together with requiring (\ref{3.16}) and (\ref{3.18}) simultaneously, implies that
\bea &&\left|m_i\bb_i+\sum_{j=1}^{\mu_i}E(m_i\aa_{i,j})-\varrho_i N-\Delta_i\right|\nn\\
          &&\qquad\quad\le \left|\left\{\frac{\varrho_i N}{M D_k}\right\}-\chi_i\right|M|D_i|+\mu_i\delta,\quad\forall\ 1\le i\le q.\lb{3.19}\eea

Notice that $\dl\max\limits_{1\le i\le q}\mu_i<\frac{1}{2}$ holds by assumption. So by (\ref{3.19}), in order to obtain (\ref{3.12})
we need to choose $M,N\in\N$ and $\chi_i$s such that the following estimate holds
\bea \left|\left\{\frac{\varrho_i N}{M D_i}\right\}-\chi_i\right|M|D_i|<\frac{1}{2}.\lb{3.20}\eea

On the other hand, by the choice (\ref{3.17}) of $m_i$, we have
\bea \{m_i\alpha_{i,j}\}&=&\left\{\left(\left[\frac{\varrho_i N}{M D_i}\right]+\chi_i\right)M\alpha_{i,j}\right\}\nn\\
&=&\left\{\frac{\varrho_i N\alpha_{i,j}}{D_i}+\left(\chi_i-\left\{\frac{\varrho_i N}{M D_i}\right\}\right)M\alpha_{i,j}\right\}\nn\\
&\equiv&\{A_{i,j}(\varrho_i N)+B_{i,j}(\varrho_i N)\}, j=1,\cdots,\mu_i, 1\le i\le q,\lb{3.21}\eea
where
\bea A_{i,j}(\varrho_i N)=\left\{\frac{\varrho_i N\alpha_{i,j}}{D_i}\right\}-\chi_{i,j},
B_{i,j}(\varrho_i N)=\left(\chi_{i}-\left\{\frac{\varrho_i N}{M D_i}\right\}\right)M\alpha_{i,j},\lb{3.22}\eea
and $\chi_{i,j}\in\{0,1\}$ will be determined later.

Following the arguments between (4.18) and (4.20) of \cite{LoZ}, it can be easily shown that $\{m_i\alpha_{i,j}\}$ must be close
enough to $0$ or $1$, i.e., satisfying (\ref{3.13}), if
\bea \max\left\{|A_{i,j}(\varrho_i N)|,\ |B_{i,j}(\varrho_i N)|\right\}<\frac{\delta_1}{3},\qquad {\rm for}\ 0<\delta_1<\delta.\lb{3.23}\eea

By (\ref{3.20}) and (\ref{3.23}), in order to get (\ref{3.12})-(\ref{3.14}) we only need to choose integers
$\chi_i,\chi_{i,j}\in\{0,1\}$ and infinitely many integers $N\in\N$ such that all the quantities
\bea \left|\left\{\frac{\varrho_i N\alpha_{i,j}}{D_i}\right\}-\chi_{i,j}\right|,\qquad \left|\left\{\frac{\varrho_i N}{M D_i}\right\}-\chi_i\right| \lb{3.24}\eea
can be made simultaneously to be small enough, which can be reduced to a dynamical problem on a torus (cf. pages 233-234
of \cite{Lon4}). Here we omit rest of details in \cite{Lon4}. \hfill\hb

\vspace{2mm}

In 2002, Y. Long and C. Zhu \cite{LoZ} has established the common index jump theorem for symplectic paths, which has become one
of the main tools to study the periodic orbit problem in Hamiltonian and symplectic dynamics. In \cite{DLW} of 2016, H. Duan,
Y. Long and W. Wang further improved this theorem to an enhanced version which gives more precise index properties of
$\ga_k^{2m_k}$ and $\ga_k^{2m_k\pm m}$ with $1\le m \le \bar{m}$ for any fixed $\bar{m}$. Under the help of Theorem 3.4, following
the proofs of Theorem 3.5 in \cite{DLW}, next we further generalize this theorem to the case of admitting the existence of
symplectic paths with negative mean indices.

\medskip

{\bf Theorem 3.6.} ({\bf Generalized common index jump theorem for symplectic paths})
{\it Let $\gamma_i\in\mathcal{P}_{\tau_i}(2n)$ for $i=1,\cdots,q$ be a finite collection of symplectic paths with nonzero mean
indices $\hat{i}(\ga_i,1)$. Let $M_i=\ga_i(\tau_i)$. We extend $\ga_i$ to $[0,+\infty)$ by (\ref{3.9}) inductively.

Then for any fixed $\bar{m}\in \N$, there exist infinitely many $(q+1)$-tuples
$(N, m_1,\cdots,m_q) \in \N^{q+1}$ such that the following hold for all $1\le i\le q$ and $1\le m\le \bar{m}$,
\bea
\nu(\ga_i,2m_i-m) &=& \nu(\ga_i,2m_i+m) = \nu(\ga_i, m),   \lb{3.25}\\
i(\ga_i,2m_i+m) &=& 2\varrho_i N+i(\ga_i,m),                         \lb{3.26}\\
i(\ga_i,2m_i-m) &=&  2\varrho_i N-i(\ga_i,m)-2(S^+_{M_i}(1)+Q_i(m)),  \lb{3.27}\\
i(\ga_i, 2m_i)&=& 2\varrho_i N -(S^+_{M_i}(1)+C(M_i)-2\Delta_i),     \lb{3.28}\eea
where \bea &&\varrho_i=\left\{\begin{array}{cc}1, &{\rm if}\ \hat{i}(\ga_i,1)>0, \cr
                                     -1, &{\rm if}\  \hat{i}(\ga_i,1)<0, \end{array}\right.\qquad
\Delta_i = \sum_{0<\{m_i\th/\pi\}<\delta}S^-_{M_i}(e^{\sqrt{-1}\th}),\nn\\
&&\ Q_i(m) = \sum_{e^{\sqrt{-1}\th}\in\sg(M_i),\atop \{\frac{m_i\th}{\pi}\}
                   = \{\frac{m\th}{2\pi}\}=0}S^-_{M_i}(e^{\sqrt{-1}\th}). \lb{3.29}\eea

More precisely, by (\ref{3.17}) and (4.40), (4.41) in \cite{LoZ} , we have
\bea m_i=\left(\left[\frac{N}{M|\hat i(\gamma_i, 1)|}\right]+\chi_i\right)M,\quad\forall\  1\le i\le q,\lb{3.30}\eea
where $\chi_i=0$ or $1$ for $1\le i\le q$ and $\frac{M\theta}{\pi}\in\Z$
whenever $e^{\sqrt{-1}\theta}\in\sigma(M_i)$ and $\frac{\theta}{\pi}\in\Q$
for some $1\le i\le q$.  Furthermore, by (\ref{3.24}),
for any $\epsilon>0$, we can choose $N$ and $\{\chi_i\}_{1\le i\le q}$ such that}
\bea \left|\left\{\frac{ N}{M|\hat i(\gamma_i, 1)|}\right\}-\chi_i\right|<\epsilon,\quad\forall\  1\le i\le q.\lb{3.31}\eea

\medskip

{\bf Proof.} For $1\le i\le q$, let $\mu_i=\sum_{0<\th<2\pi}S_{M_i}^-(e^{\sqrt{-1}\th})$, $\alpha_{i,j}=\frac{\th_j}{\pi}$
where $e^{\sqrt{-1}\th_j}\in\sigma(M_i)$ for $1\le j\le\mu_i$, and
$D_i=i(\ga_i,1) + S^+_{M_i}(1) - C(M_i) + \sum_{\th\in(0,2\pi)}\frac{\th}{\pi}S^-_{M_i}(e^{\sqrt{-1}\th})$. Then Theorem
3.6 can be proved by Theorem 3.4 and using $\varrho_kN$ and $m_i=\left(\left[\frac{\varrho_i N}{M D_i}\right]+\chi_i\right)M$
to replace the corresponding $N$ and $m_i$s in the proof of Theorem 3.5 of \cite{DLW}. Here we omit all details.\hfill\hb

\medskip

{\bf Remark 3.7.} Let $l=q+\sum_{k=1}^q \mu_k$, and
\bea &&v=\left(\frac{1}{M|\hat{i}(\gamma_1,1)|},\cdots,\frac{1}{M|\hat{i}(\gamma_1,1)|},\right.\nn\\
&&\qquad\ \quad\left.\frac{\aa_{1,1}}{|\hat{i}(\gamma_1,1)|},\cdots,
\frac{\aa_{1,\mu_1}}{|\hat{i}(\gamma_1,1)|},\cdots,\frac{\aa_{q,1}}{|\hat{i}(\gamma_q,1)|},\cdots,
\frac{\aa_{q,\mu_q}}{|\hat{i}(\gamma_q,1)|}\right)\in\R^l.\nn\eea
Theorem 3.6 also shows that for any given small $\epsilon>0$ one can find a vertex
\bes \chi=(\chi_1,\cdots,\chi_q,\chi_{1,1},\cdots,\chi_{1,\mu_1},\cdots,\chi_{q,1},\cdots,\chi_{q,\mu_q})\ees
of the cube $[0,1]^l$ and infinitely many $N\in\N$ such that $|\{Nv\}-\chi|<\epsilon$ (also see the Step 2 in the proof of Theorem 4.1 on pages 346-347 of \cite{LoZ} for the existence of $\chi$).

\medskip

{\bf Theorem 3.8.} (cf. Theorem 2.1 of \cite{HuL} and Theorem 6.1 of \cite{LLo2}) {\it Suppose $\Sg\in \H_{st}(2n)$ and
$(\tau,y)\in \T(\Sigma)$. Then we have
\be i(y^m)\equiv i(m\tau,y)=i(y, m)-n,\quad \nu(y^m)\equiv\nu(m\tau, y)=\nu(y, m),
        \lb{3.32}\ee
where $m\in\N$, $i(y^m)$ and $\nu(y^m)$ are the index and nullity of $(m\tau,y)$ defined in Section 2, $i(y, m)$ and $\nu(y, m)$
are the Maslov-type index and nullity of $(m\tau,y)$ (cf. Section 5.4 of \cite{Lon3}). In particular, we have
$\hat{i}(\tau,y)=\hat{i}(y,1)$, where $\hat{i}(\tau ,y)$ is given in Section 2, $\hat{i}(y,1)$
is the mean Maslov-type index (cf. Definition 8.1 of \cite{Lon4}). Hence we denote it simply by $\hat{i}(y)$.}

\setcounter{figure}{0}
\setcounter{equation}{0}
\section{Proofs of Theorem 1.2 and Theorem 1.3}

In order to prove Theorem 1.2, let $\Sigma\in \mathcal{H}_{st}(2n)$ be a non-degenerate perfect star-shaped hypersurface
which possesses only finitely many prime closed characteristics $\{(\tau_k,y_k)\}_{k=1}^q$ with $\hat{i}(y_k,1)\neq 0$.
Note that there exist at least one closed characteristic on $\Sg$ with positive mean index by the first identity of
(\ref{2.11}) in Theorem 2.3. So without loss of generality, the following mixed mean index condition holds:

\medskip

{\bf (MMI)} {\it There exists an integer $q_0\in[1,q]$ such that $\hat{i}(y_k,1)>0$ for $1\leq k\leq q_0$ and $\hat{i}(y_k,1)<0$ for $q_0+1\leq k\leq q$.}

\medskip

Denote by $\ga_k\equiv \ga_{y_k}$ the associated symplectic path of $(\tau_k,y_k)$ for $1\le k\le q$. Then by Lemma 3.3 of \cite{HuL} and Lemma 3.2 of \cite{Lon1},
there exists $P_k\in Sp(2n)$ and $U_k\in Sp(2n-2)$ such that
\bea M_k\equiv\ga_k(\tau_k)=P_k^{-1}(N_1(1,1)\dm U_k)P_k,\qquad\forall\  1\le k\le q,\lb{4.1}\eea
where, because $\Sigma$ is non-degenerate, every $U_k$ is isotopic in $\Omega^0(U_k)$ to a matrix of the following form by Theorem 3.2 (cf. Theorem 4.7 of \cite{Wan2})
\bea
&& R(\th_1)\,\dm\,\cdots\,\dm\,R(\th_r)\,\dm\,D(\pm 2)^{\dm s} \nn\\
&&\quad \dm\,N_2(e^{\aa_{1}\sqrt{-1}},A_{1})\,\dm\,\cdots\,\dm\,N_2(e^{\aa_{r_{\ast}}\sqrt{-1}},A_{r_{\ast}})\nn\\
&&\quad \dm\,N_2(e^{\bb_{1}\sqrt{-1}},B_{1})\,\dm\,\cdots\,\dm\,N_2(e^{\bb_{r_{0}}\sqrt{-1}},B_{r_{0}}), \nn\eea
where $\frac{\th_{j}}{2\pi}\in[0,1]\setminus\Q$ for $1\le j\le r$; $\frac{\aa_{j}}{2\pi}\in[0,1]\setminus\Q$ for $1\le j\le r_{\ast}$;
$\frac{\bb_{j}}{2\pi}\in[0,1]\setminus\Q$ for $1\le j\le r_0$ and
\be r+ s +2r_{\ast} + 2r_0 = n-1. \lb{4.2}\ee

{\bf Proof of Theorem 1.2.}

\medskip

We prove Theorem 1.2 in two cases:

\medskip

{\bf Case 1.} {\it $n$ is even.}

\medskip

We continue the proof in three steps.

{\bf Step 1.} {\it The first set of iterates for the choice of the vertex $\chi$ in the cube $[0,1]^l$.}

By (MMI), we have $\hat{i}(y_k)=\hat{i}(y_k,1)>0$ for $1\le k\le q_0$ and $\hat{i}(y_k)=\hat{i}(y_k,1)<0$ for $q_0+1\le k\le q$,
which implies that $i(y_k,m)\rightarrow +\infty$ for $1\le k\le q_0$ and $i(y_k,m)\rightarrow -\infty$ for $q_0+1\le k\le q$ as $m\rightarrow +\infty$.
So the positive integer
$\bar{m}$ defined by
\bea &&\bar{m}_1=\max_{1\le k\le q_0}\{\min\{m_0\in\N\ |\ i(y_k,m+l)\ge i(y_k,l)+n+1, \nn\\
    &&\qquad\qquad\qquad\qquad\qquad\qquad\qquad\qquad\qquad\qquad\forall\ l\ge 1, m\geq m_0\}\}\nn\\
    &&\bar{m}_2=\max_{q_0+1\le k\le q}\{\min\{m_0\in\N\ |\ i(y_k,m+l)\leq i(y_k,l)-n-1, \nn\\
    &&\qquad\qquad\qquad\qquad\qquad\qquad\qquad\qquad\qquad\qquad\forall\ l\ge 1, m\geq m_0\}\}\nn\\
    &&\bar{m}=\max\{\bar{m}_1, \bar{m}_2\}\lb{4.3}\eea
is well-defined and finite.

For the integer $\bar{m}$ defined in (\ref{4.3}), it follows from Theorem 3.6 and Remark 3.7
that there exist a vertex $\chi$ of $[0,1]^l$ and infinitely many $(q+1)$-tuples $(N, m_1, \cdots, m_q)\in\N^{q+1}$ such that for any
$1\le k\le q$, there holds
\bea
\bar{m}+2 &\le& \min\{2m_k,\ 1\le k\le q\},\lb{4.4}\\
i(y_k,{2m_k-m}) &=& 2\varrho_k N-2-i(y_k,m),\quad\forall\  1\le m\le\bar{m}, \lb{4.5}\\
i(y_k,{2m_k}) &=& 2\varrho_k N-1-C(M_k)+2\Delta_k,\lb{4.6}\\
i(y_k,{2m_k+m})&=& 2\varrho_k N+i(y_k,m),\quad \forall\ 1\le m\le\bar{m},\lb{4.7}\eea
where note that $S^+_{M_k}(1)=1,Q_k(m)=0$, $\forall\ m\ge 1$ by (\ref{4.1})-(\ref{4.2}) and Lemma 3.1.

By the definition (\ref{4.3}) of $\bar{m}$ and (\ref{4.6}), for any $m\ge \bar{m}+1$, we obtain
\bea
 i(y_k,2m_k-m)&\le& i(y_k,2m_k)-n-1\nn\\&=&2N-n-2+2\Delta_k-C(M_k)\nn\\&\le& 2N-3,\quad 1\le k\le q_0, \lb{4.8}\\
 i(y_k,2m_k-m)&\ge& i(y_k,2m_k)+n+1\nn\\&=&-2N+n+2\Delta_k-C(M_k)\nn\\&\ge& -2N+1,\quad q_0+1\le k\le q, \lb{4.8'}\\
 i(y_k,2m_k+m)&\ge& i(y_k,2m_k)+n+1\nn\\&=&2N+n-C(M_k)+2\Delta_k\nn\\ &\ge& 2N+1,\quad 1\le k\le q_0, \lb{4.9}\\
  i(y_k,2m_k+m)&\le& i(y_k,2m_k)-n-1\nn\\&=&-2N-n-2-C(M_k)+2\Delta_k\nn\\ &\le& -2N-3,\quad q_0+1\le k\le q, \lb{4.9'}\eea
where we have used the fact $2\Delta_k-C(M_k)\le \Delta_k\le C(M_k)\le n-1$, which follows from $\Delta_k\le C(M_k)\le n-1$ by Lemma 3.1, the definition of $\Delta_k$ in (3.29) and the definition of $C(M_k)$ in Theorem 3.3.

Then by (\ref{4.5})-(\ref{4.9'}) and Theorem 3.8, we obtain
\bea  i(y_k^{2m_k-m})&\le& 2N-n-3,\nn\\
&&\ \forall\ \bar{m}+1\le m\le 2m_k-1,\ 1\le k\le q_0,\lb{4.10} \\
 i(y_k^{2m_k-m}) &\ge& -2N-n+1,\nn\\
 &&\ \forall\ \bar{m}+1\le m\le 2m_k-1,\ q_0+1\le k\le q,\lb{4.10'} \\
 i(y_k^{2m_k-m}) &=& 2\varrho_k N-2n-2-i(y_k^m),\quad \forall\ 1\le m\le \bar{m}, \lb{4.11}\\
 i(y_k^{2m_k}) &=& 2\varrho_k N-C(M_k)+2\Delta_k-n-1,\lb{4.12}\\
 i(y_k^{2m_k+m}) &=& 2\varrho_k N+i(y_k^m),\quad \forall\ 1\le m\le \bar{m},\lb{4.13}\\
 i(y_k^{2m_k+m}) &\ge& 2N-n+1,\quad \forall\ m\ge \bar{m}+1,\ 1\le k\le q_0,\lb{4.14} \\
i(y_k^{2m_k+m})&\le& -2N-n-3,\quad \forall\ m\ge \bar{m}+1,\ q_0+1\le k\le q.\lb{4.14'}\eea

{\bf Claim 1:} {\it One can choose $N\in\N$ in Theorem 3.6 satisfying (\ref{4.10})-(\ref{4.14'}) such that
\bea \sum_{k=1}^q 2m_k\hat{\chi}(y_k)=N.\lb{4.15}\eea}

In fact, let $\ep<\frac{1}{1+2M\sum_{1\le k\le q}|\hat{\chi}(y_k)|}$, by Theorem 2.3 and (MMI) we have
\bea \sum_{k=1}^q\frac{\hat{\chi}(y_k)}{|\hat{i}(y_k)|}=\sum_{\hat{i}(y_k)>0}\frac{\hat{\chi}(y_k)}{\hat{i}(y_k)}-
\sum_{\hat{i}(y_k)<0}\frac{\hat{\chi}(y_k)}{\hat{i}(y_k)}=\frac{1}{2},\nn\eea
which, together with (\ref{3.30})-(\ref{3.31}) and recalling that $\hat{i}(y_k,1)=\hat{i}(y_k)$, yields
\bea \left|N-\sum_{k=1}^q 2m_k\hat{\chi}(y_k)\right|
&=& \left|\sum_{k=1}^q\frac{2N\hat{\chi}(y_k)}{|\hat{i}(y_k)|}-\sum_{k=1}^q 2\hat{\chi}(y_k)
    \left(\left[\frac{N}{M|\hat{i}(y_k)|}\right]+\chi_k\right)M\right|\nn\\
&\le& 2M\sum_{k=1}^q |\hat{\chi}(y_k)|\left|\left\{\frac{N}{M|\hat{i}(y_k)|}\right\}-\chi_k\right|.\nn\\
&<& 2M\ep\sum_{k=1}^q|\hat{\chi}(y_k)|\nn\\
&<& 1.\lb{4.16}\eea
It proves Claim 1 since each $2m_k\hat{\chi}(y_k)$ is an integer.

\medskip

Now by Lemma 2.2, Theorem 3.3 and Theorem 3.8, it yields
\bea
&& \sum_{m=1}^{2m_k} (-1)^{d(K)+i(y_k^m)}\dim C_{S^1,d(K)+i(y_k^{m})}(F_{a, K},S^1\cdot x_k^m)\nn\\
&&\quad = \sum_{m=1}^{2m_k} (-1)^{i(y_k^m)}\dim C_{S^1,d(K)+i(y_k^{m})}(F_{a, K},S^1\cdot x_k^m)\nn\\
&&\quad= \sum_{j=0}^{m_k-1} \sum_{m=2j+1}^{2j+2} (-1)^{i(y_k^m)} \dim C_{S^1,d(K)+i(y_k^{m})}(F_{a, K},S^1\cdot x_k^m)\nn\\
&&\quad = \sum_{j=0}^{m_k-1} \sum_{m=1}^{2} (-1)^{i(y_k^m)} \dim C_{S^1,d(K)+i(y_k^{m})}(F_{a, K},S^1\cdot x_k^m)\nn\\
&&\quad = m_k \sum_{m=1}^{2} (-1)^{i(y_k^m)} \dim C_{S^1,d(K)+i(y_k^{m})}(F_{a, K},S^1\cdot x_k^m)\nn\\
&&\quad = 2m_k\hat{\chi}(y_k),\qquad \forall\ 1\le k\le q, \lb{4.17}\eea
where $x_k$ is the critical point of $F_{a, K}$ corresponding to $y_k$, and we choose large enough $K$ such that
$d(K)=2n([KT/{2\pi}]+1)\ge -i(y_k^{m})$ for $1\le m\le 2m_k$ and $1\le k\le q$. In addition, the third equality follows from Lemma 2.2 and the fact $i(y_k^{m+2})-i(y_k^m)\in 2\Z$ for any $m\in\N$ from (\ref{3.10}) of Theorem 3.3 and Theorem 3.8.

For $1\le k\le q$, by (\ref{4.14})-(\ref{4.14'}) and Lemma 2.2, we know that all $y_k^{2m_k+m}$'s with $m\ge\bar{m}+1$
have no contribution to the alternative sum $\sum_{p=-2N-n-1}^{2N-n-1}(-1)^p M_p$, where the Morse-type number $M_p$
is defined in (\ref{2.14}). Similarly again by
Lemma 2.2 and (\ref{4.10})-(\ref{4.10'}), all $y_k^{2m_k-m}$'s with $\bar{m}+1\le m\le 2m_k-1$ only have contribution to
$\sum_{p=-2N-n-1}^{2N-n-1}(-1)^p M_p$.

For $1\le m\le\bar{m}$, by (\ref{4.13}) and Lemma 2.2, we know that all $y_k^{2m_k+m}$'s with $-n\le i(y_k^m)$
for $1\le k\le q_0$, or $ i(y_k^m)\leq -n-2$
for $q_0+1\le k\le q$, have no contribution to the alternative sum $\sum_{p=-2N-n-1}^{2N-n-1}(-1)^p M_p$.
Similarly again by Lemma 2.2 and (\ref{4.11}), for $1\le m\le\bar{m}$, all $y_k^{2m_k-m}$'s with $-n\le i(y_k^m)$ and
$1\le k\le q_0$, or $ i(y_k^m)\leq -n-2$
for $q_0+1\le k\le q$, only have contribution to $\sum_{p=-2N-n-1}^{2N-n-1}(-1)^p M_p$.

Since $i(y_k^{m})\neq -n-1$ when $(m\tau_k, y_k)$ is good, by (MMI), Definition 1.1 and Theorem 3.8, we set
\bea
M_+^e(k)=\left\{\begin{array}{ll}^{\#}\{1\le m\le \bar{m}\ |\ i(y_k^{m})\le -n-2,\nn\\
 \qquad\  i(y_k^{2m_k+m})\in 2\Z, i(y_k)\in 2\Z\}, ~{\rm if}~ 1\le k\le q_0, \\
^{\#}\{1\le m\le \bar{m}\ |\ i(y_k^{m})\geq -n,\nn\\
 \qquad\  i(y_k^{2m_k+m})\in 2\Z, i(y_k)\in 2\Z\},~{\rm if}~q_0+1\le k\le q, \end{array}\right. \nn
\eea
\bea
M_+^o(k)=\left\{\begin{array}{ll}^{\#}\{1\le m\le \bar{m}\ |\ i(y_k^{m})\le -n-2,\nn\\
  \qquad\  i(y_k^{2m_k+m})\in 2\Z-1,\ i(y_k)\in 2\Z-1\}, ~{\rm if}~ 1\le k\le q_0, \\ ^{\#}\{1\le m\le \bar{m}\ |\ i(y_k^{m})\geq -n,\ i(y_k^{2m_k+m})\in 2\Z-1,\nn\\
   \qquad\  i(y_k)\in 2\Z-1\},~{\rm if}~ q_0+1\le k\le q, \end{array}\right. \nn\eea
\bea
M_-^e(k)=\left\{\begin{array}{ll}^{\#}\{1\le m\le \bar{m}\ |\ i(y_k^{m})\le -n-2,\nn\\
 \qquad\  i(y_k^{2m_k-m})\in 2\Z,\ i(y_k)\in 2\Z\}, ~{\rm if}~ 1\le k\le q_0,
 \\ ^{\#}\{1\le m\le \bar{m}\ |\ i(y_k^{m})\geq -n,\nn\\
 \qquad\  i(y_k^{2m_k-m})\in 2\Z,\ i(y_k)\in 2\Z\},~{\rm if}~q_0+1\le k\le q, \end{array}\right. \nn\eea
\bea
M_-^o(k)=\left\{\begin{array}{ll}^{\#}\{1\le m\le \bar{m}\ |\ i(y_k^{m})\le -n-2,\nn\\
 \quad\  i(y_k^{2m_k-m})\in 2\Z-1,\ i(y_k)\in 2\Z-1\}, ~{\rm if}~ 1\le k\le q_0,
 \\  ^{\#}\{1\le m\le \bar{m}\ |\ i(y_k^{m})\geq -n,\nn\\
  \quad\  i(y_k^{2m_k-m})\in 2\Z-1,\ i(y_k)\in 2\Z-1\},~{\rm if}~q_0+1\le k\le q, \end{array} \right. \nn\eea
which, together with $i(y_k^{2m_k+m})-i(y_k^{2m_k-m})\in 2\Z$ by (\ref{4.11}) and (\ref{4.13}), yields
\be M_+^{e}(k)=M_-^{e}(k),\quad M_+^{o}(k)=M_-^{o}(k),\qquad\forall\  1\le k\le q.\lb{4.18}\ee

Here by Lemma 2.2 we remark that $M_\pm^e(k)$ with $1\le k\le q$ only counts the contributions of $y_k^{2m_k\pm m}$ with $1\le m\le \bar{m}$ to the even-th Morse-type numbers, and $M_\pm^o(k)$ only counts the contributions of $y_k^{2m_k\pm m}$ with $1\le m\le \bar{m}$ to the odd-th Morse-type numbers. Next we focus on these precise contributions to the alternative sum $\sum_{p=-2N-n-1}^{2N-n-1}(-1)^p M_p$.

For $1\le k\le q_0$ and $1\le m\le \bar{m}$ satisfying $i(y_k^{m})\le -n-2$, and for
$q_0+1\le k\le q$ and $1\le m\le \bar{m}$ satisfying $i(y_k^{m})\geq -n$, by (\ref{4.11}) and (\ref{4.13})
it yields
\bea i(y_k^{2m_k-m})&\ge& 2N-n,\quad i(y_k^{2m_k+m})\le 2N-n-2, \nn\\
&&\qquad\qquad\qquad\qquad \forall\ 1\le k\le q_0,\lb{4.19}\\
i(y_k^{2m_k-m})&\le& -2N-n-2,\quad i(y_k^{2m_k+m})\ge -2N-n, \nn\\
&&\qquad\qquad\qquad\qquad \forall\  q_0+1\le k\le q.\lb{4.19'}\eea

So, for $1\le m\le\bar{m}$, by (\ref{4.19})-(\ref{4.19'}) and Lemma 2.2, we know that all $y_k^{2m_k+m}$'s (i.e., $M_+^{e,o}(k)$) with $i(y_k^m)\le -n-2$ and $1\le k\le q_0$, or $i(y_k^{m})\geq -n$ and $q_0+1\le k\le q$, only have contribution to the alternative
sum $\sum_{p=-2N-n-1}^{2N-n-1}(-1)^p M_p$, and all $y_k^{2m_k-m}$'s (i.e., $M_-^{e,o}(k)$) with $i(y_k^m)\le -n-2$ and $1\le k\le q_0$, or  $i(y_k^{m})\geq -n$ and $q_0+1\le k\le q$, have no contribution to $\sum_{p=-2N-n-1}^{2N-n-1}(-1)^p M_p$.

Thus for the Morse-type numbers $M_p$'s in (\ref{2.14}), by (\ref{4.17})-(\ref{4.19'}) we have
{\small \bea &&\sum_{p=-2N-n-1}^{2N-n-1}(-1)^p M_p\nn\\
&&= \sum_{k=1}^{q}\ \sum_{1\le m\le 2m_k+\bar{m} \atop -2N-n-1\le i(y_k^{m})\le 2N-n-1}
    (-1)^{d(K)+i(y_k^m)}\dim C_{S^1,d(K)+i(y_k^{m})}(F_{a, K},S^1\cdot x_k^m)\nn\\
&&= \sum_{k=1}^{q}\ \sum_{m=1}^{2m_k} (-1)^{d(K)+i(y_k^m)}\dim C_{S^1,d(K)+i(y_k^{m})}(F_{a, K},S^1\cdot x_k^m)\nn\\
&&\quad +\sum_{k=1}^q \left[M_+^e(k)-M_+^o(k)\right]-\sum_{k=1}^q \left[M_-^e(k)-M_-^o(k)\right]\nn\\
&&\quad -\sum_{1\le k\le q_0 \atop i(y_k^{2m_k})\ge 2N-n} (-1)^{i(y_k^{2m_k})}
              \dim C_{S^1,d(K)+i(y_k^{2m_k})}(F_{a, K},S^1\cdot x_k^{2m_k})\nn\\
&&\quad -\sum_{q_0+1\le k\le q \atop i(y_k^{2m_k})\le -2N-n-2} (-1)^{i(y_k^{2m_k})}
              \dim C_{S^1,d(K)+i(y_k^{2m_k})}(F_{a, K},S^1\cdot x_k^{2m_k})\nn\\
&&= \sum_{k=1}^{q} 2m_k\hat{\chi}(y_k)\nn\\
&&\quad -\sum_{1\le k\le q_0 \atop i(y_k^{2m_k})\ge 2N-n} (-1)^{i(y_k^{2m_k})}
              \dim C_{S^1,d(K)+i(y_k^{2m_k})}(F_{a, K},S^1\cdot x_k^{2m_k})\nn\\
&&\quad -\sum_{q_0+1\le k\le q \atop i(y_k^{2m_k})\le -2N-n-2} (-1)^{i(y_k^{2m_k})}
              \dim C_{S^1,d(K)+i(y_k^{2m_k})}(F_{a, K},S^1\cdot x_k^{2m_k}). \lb{4.20}\eea}

In order to exactly know whether the iterate $y_k^{2m_k}$ of $y_k$ has contribution to the alternative sum
$\sum_{p=-2N-n-1}^{2N-n-1}(-1)^p M_p$,  $1\le k\le q$, we set
\bea
&&N_+^e =~^{\#}\{q_0+1\le k\le q\ |\ i(y_k^{2m_k})\le -2N-n-2,\nn\\
&&\qquad\qquad\qquad\qquad i(y_k^{2m_k})\in 2\Z,\ i(y_k)\in 2\Z\}\nn\\
&&\qquad\quad+ ~^{\#}\{1\le k\le q_0\ |\ i(y_k^{2m_k})\ge 2N-n,\nn\\
&&\qquad\qquad\qquad\qquad i(y_k^{2m_k})\in 2\Z,\ i(y_k)\in 2\Z\},\lb{4.21}
\eea
\bea
&&N_+^o =~^{\#}\{q_0+1\le k\le q\ |\ i(y_k^{2m_k})\le -2N-n-2,\nn\\
&&\qquad\qquad\qquad\qquad i(y_k^{2m_k})\in 2\Z-1,\ i(y_k)\in 2\Z-1\}\nn\\
&&\qquad\quad +~ ^{\#}\{1\le k\le q_0\ |\ i(y_k^{2m_k})\ge 2N-n,\nn\\
&&\qquad\qquad\qquad\qquad i(y_k^{2m_k})\in 2\Z-1,\ i(y_k)\in 2\Z-1\},\lb{4.22}
\eea
\bea
&&N_-^e =~^{\#}\{q_0+1\le k\le q\ |\ i(y_k^{2m_k})\ge -2N-n,\nn\\
&&\qquad\qquad\qquad\qquad i(y_k^{2m_k})\in 2\Z,\ i(y_k)\in 2\Z\}\nn\\
&&\qquad\quad +~^{\#}\{1\le k\le q_0\ |\ i(y_k^{2m_k})\le 2N-n-2,\nn\\
&&\qquad\qquad\qquad\qquad i(y_k^{2m_k})\in 2\Z,\ i(y_k)\in 2\Z\},\lb{4.23}
\eea
\bea
&&N_-^o =~^{\#}\{q_0+1\le k\le q\ |\ i(y_k^{2m_k})\ge -2N-n,\nn\\
&&\qquad\qquad\qquad\qquad i(y_k^{2m_k})\in 2\Z-1,\ i(y_k)\in 2\Z-1\}\nn\\
&&\qquad\quad +~^{\#}\{1\le k\le q_0\ |\ i(y_k^{2m_k})\le 2N-n-2,\nn\\
&&\qquad\qquad\qquad\qquad i(y_k^{2m_k})\in 2\Z-1,\ i(y_k)\in 2\Z-1\}.\lb{4.24}\eea

Thus by Claim 1, (\ref{4.20}), the definitions of $N^{e}_+$ and $N_+^o$ and (\ref{2.16}), we have
\bea N+N_+^o-N_+^e
&=& \sum_{k=1}^q 2m_k\hat{\chi}(y_k)+N_+^o-N_+^e  \nn\\
&=& \sum_{p=-2N-n-1}^{2N-n-1}(-1)^p M_p \nn\\
&\le& \sum_{p=-2N-n-1}^{2N-n-1}(-1)^p b_p=\sum_{p=0}^{2N-n-2}b_p \nn\\
&=& N-\frac{n}{2},\lb{4.25}\eea
where the first equality holds by Claim 1, the second equality follows from (4.25) and the definitions of $N^{e}_+$
and $N_+^o$, and the last equality follows from $b_{2j}=1$ and $b_{2j-1}=0$ for $0\le j\le N-\frac{n-2}{2}$ by
(\ref{2.16}) where $n$ is even.

So (\ref{4.25}) give the following estimate
\be N_+^e\ge \frac{n}{2}.\lb{4.26}\ee

{\bf Step 2.} {\it The second set of iterates for the choice of the dual vertex $\hat{\chi}=1-\chi$ in the cube $[0,1]^l$.}

Similar to (\ref{4.10})-(\ref{4.14'}), for $\hat{\chi}=\hat{1}-\chi$ of the cube $[0,1]^l$ with $\chi$ chosen below
(\ref{4.3}) where $\hat{1}=(1,\cdots,1)$, it follows from Theorem 3.6 (also cf. Theorem 2.8 of \cite{HaW} and Theorem
4.2 of \cite{LoZ}) and Remark 3.7 that there exist also infinitely many $(q+1)$-tuples
$(N', m_1', \cdots, m_q')\in\N^{q+1}$ such that for any $1\le k\le q$, there holds
\bea
i(y_k^{2m_k'-m})&\le& 2N'-n-3,\nn\\
&&\ \forall\ \bar{m}+1\le m\le 2m_k'-1,\ 1\le k\le q_0,\lb{4.27} \\
i(y_k^{2m_k'-m})&\ge& -2N'-n+1,\nn\\
&&\ \forall\ \bar{m}+1\le m\le 2m_k'-1,\ q_0+1\le k\le q,\lb{4.27'} \\
i(y_k^{2m_k'-m}) &=& 2\varrho_k N'-2n-2-i(y_k^m),\quad \forall\ 1\le m\le \bar{m}, \lb{4.28}\\
i(y_k^{2m_k'}) &=& 2\varrho_k N'-C(M_k)+2\Delta_k'-n-1,\lb{4.29}\\
i(y_k^{2m_k'+m}) &=& 2\varrho_k N'+i(y_k^m),\quad \forall\ 1\le m\le \bar{m},\lb{4.30}\\
i(y_k^{2m_k'+m}) &\ge& 2N'-n+1,\quad \forall\ m\ge \bar{m}+1,\ 1\le k\le q_0,\lb{4.31}\\
i(y_k^{2m_k'+m}) &\le& -2N'-n-3,\  \forall\ m\ge \bar{m}+1,\ q_0+1\le k\le q,\lb{4.31'} \eea
where, furthermore, $\Delta_k$ and $\Delta_k'$ satisfy the following relationship
\bea \Delta_k' + \Delta_k = C(M_k),\qquad\forall\ 1\le k\le q,\lb{4.32}\eea
by the fact $\hat{\chi}=\hat{1}-\chi$ and the proof of Claim 4 in the proof of Theorem 1.1 of \cite{DLW} or (42) in Theorem 2.8 of \cite{HaW}.

Similarly, we define
\bea
&&N_+^{'e} =~^{\#}\{q_0+1\le k\le q\ |\ i(y_k^{2m_k'})\le -2N'-n-2,\nn\\
&&\qquad\qquad\qquad\qquad i(y_k^{2m_k'})\in 2\Z,\ i(y_k)\in 2\Z\}\nn\\
&&\qquad\quad+~^{\#}\{1\le k\le q_0\ |\ i(y_k^{2m_k'})\ge 2N'-n,\nn\\
&&\qquad\qquad\qquad\qquad i(y_k^{2m_k'})\in 2\Z,\ i(y_k)\in 2\Z\},   \lb{4.33}
\eea
\bea
&&N_+^{'o} =~^{\#}\{q_0+1\le k\le q\ |\ i(y_k^{2m_k'})\le -2N'-n-2,\nn\\
&&\qquad\qquad\qquad\qquad i(y_k^{2m_k'})\in 2\Z-1,\ i(y_k)\in 2\Z-1\}\nn\\
&&\qquad\quad+~^{\#}\{1\le k\le q_0\ |\ i(y_k^{2m_k'})\ge 2N'-n,\nn\\
&&\qquad\qquad\qquad\qquad i(y_k^{2m_k'})\in 2\Z-1,\ i(y_k)\in 2\Z-1\}, \lb{4.34}
\eea
\bea
&&N_-^{'e} =~^{\#}\{q_0+1\le k\le q\ |\ i(y_k^{2m_k'})\ge -2N'-n,\nn\\
&&\qquad\qquad\qquad\qquad i(y_k^{2m_k'})\in 2\Z,\ i(y_k)\in 2\Z\}\nn\\
 &&\qquad\quad+~^{\#}\{1\le k\le q_0\ |\ i(y_k^{2m_k'})\le 2N'-n-2,\nn\\
&&\qquad\qquad\qquad\qquad i(y_k^{2m_k'})\in 2\Z,\ i(y_k)\in 2\Z\},   \lb{4.35}
\eea
\bea
&&N_-^{'o} =~^{\#}\{q_0+1\le k\le q\ |\ i(y_k^{2m_k'})\ge -2N'-n,\nn\\
&&\qquad\qquad\qquad\qquad i(y_k^{2m_k'})\in 2\Z-1,\ i(y_k)\in 2\Z-1\}\nn\\
&&\qquad\quad+~^{\#}\{1\le k\le q_0\ |\ i(y_k^{2m_k'})\le 2N'-n-2,\nn\\
&&\qquad\qquad\qquad\qquad i(y_k^{2m_k'})\in 2\Z-1,\ i(y_k)\in 2\Z-1\}. \lb{4.36}
\eea
So by (\ref{4.29}) and (\ref{4.32}) it yields
\bea i(y_k^{2m_k'}) &=& 2\varrho_k N'-C(M_k)+2(C(M_k)-\Delta_k)-n-1\nn\\
&=&2\varrho_k N'+C(M_k)-2\Delta_k-n-1. \lb{4.37}\eea

So by definitions (\ref{4.21})-(\ref{4.24}) and (\ref{4.33})-(\ref{4.36}) we have
\be N_{\pm}^e=N_{\mp}^{'e},\qquad N_{\pm}^o=N_{\mp}^{'o}. \lb{4.38}\ee

Thus, carrying out the arguments similar to (\ref{4.25})-(\ref{4.26}), by Claim 1, the definitions of
$N^{'e}_+$ and $N^{'o}_+$  and (\ref{2.16}), we have
\bea N'+N_+^{'o}-N_+^{'e}&=&\sum_{k=1}^q 2m'_k\hat{\chi}(y_k)+N_+^{'o}-N_+^{'e}\nn\\
&=&\sum_{p=-2N'-n-1}^{2N'-n-1}(-1)^p M_p\nn\\
&\le&\sum_{p=-2N'-n-1}^{2N'-n-1}(-1)^p b_p=\sum_{p=0}^{2N'-n-2}b_p\nn\\
&=& N'-\frac{n}{2},\lb{4.39}\eea
which, together with (\ref{4.38}), implies
\bea N_-^{e}=N_+^{'e}\ge\frac{n}{2}.\lb{4.40}\eea

{\bf Step 3.} {\it The summary.}

So by (\ref{4.26}) and (\ref{4.40}) it yields
\be q\ge N_+^{e}+N_-^{e}\ge n. \lb{4.41}\ee

In addition, any hyperbolic closed characteristic $y_k$ must have $i(y_k^{2m_k})=2\varrho_k N-n-1$ since there holds
$C(M_k)=0$ in the hyperbolic case. However, by (\ref{4.21}) and (\ref{4.23}), there exist at least
$(N_+^{e}+N_-^{e})$ closed characteristics with even indices $i(y_k^{2m_k})$. So all these
$(N_+^{e}+N_-^{e})$ closed characteristics are non-hyperbolic. Then (\ref{4.41}) shows that there exist
at least $n$ distinct non-hyperbolic closed characteristics. Now (\ref{4.21}), (\ref{4.23}) and (\ref{4.41})
show that $i(y_k)\in 2\Z, i(y_k^{2m_k})\in 2\Z$ for all these non-hyperbolic closed characteristics $y_k$ with $1\le k\le n$. Note that $i(y_k^{m+2})-i(y_k^m)\in 2\Z,\forall m\in\N$ by (\ref{3.10}) of Theorem 3.3 and Theorem 3.8, therefore all these $n$ distinct non-hyperbolic closed characteristics and their iterations have even Maslov-type indices. This completes the proof of Case 1.

\medskip

{\bf Case 2.} {\it $n$ is odd.}

\medskip

In this case, (MMI) still holds. Here the arguments are similar to those in the proof of Case 1. So we
only give some different parts in the proof and omit other details.

\medskip

{\bf Claim 2:} {\it There exist at least $(n-1)$ geometrically distinct non-hyperbolic closed characteristics
denoted by $\{y_k\}_{k=1}^{n-1}$ with odd Maslov-type indices on such hypersurface $\Sg$.}

\medskip

Here one crucial and different point from the proof of Case 1 is that we need to consider the alternative sum
$\sum_{p=-2N-n}^{2N-n}(-1)^p M_p$ (cf. (\ref{4.44})) instead of $\sum_{p=-2N-n-1}^{2N-n-1}(-1)^p M_p$
(cf. (\ref{4.20})). This difference is mainly due to the different parity of $n$. Since the method is similar
to that in proof of Case 1, we only list some necessary parts.

At first, there holds $i(y_k^m)\notin\{ -n-2,-n-1,-n\}$ when $(m\tau_k, y_k)$ is good, by (MMI), Definition 1.1 and Theorem 3.8, we set
\bea
\bar{M}_+^e(k)=\left\{\begin{array}{ll}^{\#}\{1\le m\le \bar{m}\ |\ i(y_k^{m})\le -n-3,\nn\\
\qquad\  i(y_k^{2m_k+m})\in 2\Z,\ i(y_k)\in 2\Z\}, ~{\rm if}~ 1\le k\le q_0, \\
              ^{\#}\{1\le m\le \bar{m}\ |\ i(y_k^{m})\geq -n+1,\nn\\
\qquad\ i(y_k^{2m_k+m})\in 2\Z,\ i(y_k)\in 2\Z\},~{\rm if}~q_0+1\le k\le q, \end{array}\right. \nn\eea
\bea
\bar{M}_+^o(k)=\left\{\begin{array}{ll}^{\#}\{1\le m\le \bar{m}\ |\ i(y_k^{m})\le -n-3,\nn\\
\qquad\ i(y_k^{2m_k+m})\in 2\Z-1,\ i(y_k)\in 2\Z-1\}, ~{\rm if}~ 1\le k\le q_0, \\
              ^{\#}\{1\le m\le \bar{m}\ |\ i(y_k^{m})\geq -n+1,\nn\\
\qquad\ i(y_k^{2m_k+m})\in 2\Z-1,\ i(y_k)\in 2\Z-1\},~{\rm if}~q_0+1\le k\le q, \end{array}\right. \nn\eea
\bea
\bar{M}_-^e(k)=\left\{\begin{array}{ll}^{\#}\{1\le m\le \bar{m}\ |\ i(y_k^{m})\le -n-3,\nn\\
\qquad\ i(y_k^{2m_k-m})\in 2\Z,\ i(y_k)\in 2\Z\}, ~{\rm if}~ 1\le k\le q_0, \\
              ^{\#}\{1\le m\le \bar{m}\ |\ i(y_k^{m})\geq -n+1,\nn\\
\qquad\ i(y_k^{2m_k-m})\in 2\Z,\ i(y_k)\in 2\Z\},~{\rm if}~q_0+1\le k\le q, \end{array}\right. \nn\eea
\bea
\bar{M}_-^o(k)=\left\{\begin{array}{ll}^{\#}\{1\le m\le \bar{m}\ |\ i(y_k^{m})\le -n-3,\nn\\
\quad\ i(y_k^{2m_k-m})\in 2\Z-1,\ i(y_k)\in 2\Z-1\}, ~{\rm if}~ 1\le k\le q_0, \\
              ^{\#}\{1\le m\le \bar{m}\ |\ i(y_k^{m})\geq -n+1,\nn\\
\quad\ i(y_k^{2m_k-m})\in 2\Z-1,\ i(y_k)\in 2\Z-1\},~{\rm if}~q_0+1\le k\le q, \end{array}\right. \nn\eea
which, together with $i(y_k^{2m_k+m})-i(y_k^{2m_k-m})\in 2\Z$ by (\ref{4.11}) and (\ref{4.13}), yields
\bea \bar{M}_+^{e}(k)=\bar{M}_-^{e}(k),\quad \bar{M}_+^{o}(k)=\bar{M}_-^{o}(k),\qquad\forall\  1\le k\le q.\lb{4.42}\eea

For $1\le k\le q_0$ and $1\le m\le \bar{m}$ satisfying $i(y_k^{m})\le -n-3$, and for
$q_0+1\le k\le q$ and $1\le m\le \bar{m}$ satisfying $i(y_k^{m})\geq -n+1$, by (\ref{4.11}) and (\ref{4.13})
it yields
\bea i(y_k^{2m_k-m})&\ge& 2N-n+1,\quad i(y_k^{2m_k+m})\le 2N-n-3,\nn\\
&&\qquad\qquad\qquad\qquad\forall\  1\le k\le q_0,\lb{4.43}\\
i(y_k^{2m_k-m})&\le& -2N-n-3,\quad i(y_k^{2m_k+m})\ge -2N-n+1, \nn\\
&&\qquad\qquad\qquad\qquad\forall \ q_0+1\le k\le q.\lb{4.43'}\eea

Then, similarly to the equation (\ref{4.20}), we have
\bea &&\sum_{p=-2N-n}^{2N-n}(-1)^p M_p = \sum_{k=1}^{q} 2m_k\hat{\chi}(y_k)\nn\\
&&\quad -\sum_{1\le k\le q_0 \atop i(y_k^{2m_k}) \ge 2N-n+1} (-1)^{i(y_k^{2m_k})}
      \dim C_{S^1,d(K)+i(y_k^{2m_k})}(F_{a, K},S^1\cdot x_k^{2m_k})\nn\\
&&\quad -\sum_{q_0+1\le k\le q \atop i(y_k^{2m_k}) \le -2N-n-3} (-1)^{i(y_k^{2m_k})}
      \dim C_{S^1,d(K)+i(y_k^{2m_k})}(F_{a, K},S^1\cdot x_k^{2m_k}).\nn\\ \lb{4.44}\eea

Denote by $H_+^e,H_+^o,H_-^e,H_-^o$ the numbers similarly defined by (\ref{4.21})-(\ref{4.24}) where $\pm2N-n$ and
$\pm2N-n-2$ are replaced by $\pm2N-n+1$ and $\pm2N-n-3$, respectively.

Then by Claim 1, (\ref{4.44}), the definitions of $H^{e}_+$ and $H_+^o$, and (\ref{2.16}), we have
\bea N+H_+^o-H_+^e&=&\sum_{k=1}^q 2m_k\hat{\chi}(y_k)+H_+^o-H_+^e\nn\\
&=&\sum_{p=-2N-n}^{2N-n}(-1)^p M_p \nn\\
&\le&\sum_{p=-2N-n}^{2N-n}(-1)^p b_p\nn\\
&=&\sum_{p=0}^{2N-n-1} b_p=\frac{2N-n-1}{2}+1\nn\\
&=& N-\frac{n-1}{2},\lb{4.45}\eea
which yields
\bea H_+^e\ge H_+^e-H_+^o\ge \frac{n-1}{2}.\lb{4.46}\eea

Similarly, denote by $H_+^{'e},H_+^{'o},H_-^{'e},H_-^{'o}$ the numbers similarly defined by
(\ref{4.33})-(\ref{4.36}) where $\pm2N'-n$ and $\pm2N'-n-2$ are replaced by $\pm2N'-n+1$ and $\pm2N'-n-3$, respectively,
and these numbers satisfy the following relationship
\be H_{\pm}^e=H_{\mp}^{'e},\qquad H_{\pm}^o=H_{\mp}^{'o}. \lb{4.47}\ee

Similarly to the inequality (\ref{4.39}), by the same arguments above and (\ref{4.47}) we can obtain
\bea H_-^e=H_+^{'e}\ge H_+^{'e}-H_+^{'o}\ge \frac{n-1}{2}.\lb{4.48}\eea

Therefore it follows from (\ref{4.46}) and (\ref{4.48}) that
\bea q\ge H_+^e+H_-^e\ge \frac{n-1}{2}+\frac{n-1}{2}=n-1.\lb{4.49}\eea

By the same arguments in the proof of Case 1, it follows from the definitions of $H_+^e$ and $H_-^e$ that
these $(n-1)$ distinct closed geodesics are non-hyperbolic, and the Viterbo indices of them and their iterations
are even, and thus the Maslov-type indices of them and their iterations are odd. This completes the proof of
Claim 2.

\medskip

{\bf Claim 3:} {\it There exist at least another geometrically distinct closed characteristic different from
those found in Claim 2 with odd Maslov-type indices on such hypersurface $\Sg$.}

\medskip

In fact, for those $(n-1)$ distinct closed characteristics $\{y_k\}_{k=1}^{n-1}$ found in Claim 2, there holds
$i(y_k^{2m_k})\neq \pm2N-n-1$ by  the definitions of $H_+^e$ and $H_-^e$, which, together with
(\ref{4.10})-(\ref{4.14'}) and (MMI), yields
\bea i(y_k^m)\neq 2N-n-1,\qquad\forall\  m\ge 1,\quad k=1,\cdots,n-1.\lb{4.50}\eea
Then by Lemma 2.2 it yields
\bea \sum_{1\le k\le n-1 \atop m\ge 1} \dim C_{S^1,d(K)+2N-n-1}(F_{a, K},S^1\cdot x_k^m)=0.\lb{4.51} \eea

By (\ref{4.10})-(\ref{4.11}), (\ref{4.13})-(\ref{4.14}) and (MMI), it yields $i(y_k^m)\neq 2N-n-1$
for any $m\neq 2m_k$ and $1\le k\le q$. Therefore, by (\ref{4.51}) and (\ref{2.14}) we obtain
\bea &&\sum_{n\le k\le q}\dim C_{S^1,d(K)+2N-n-1}(F_{a, K},S^1\cdot x_k^{2m_k})\nn\\
&&\qquad = \sum_{n\le k\le q,\ m\ge 1}\dim C_{S^1,d(K)+2N-n-1}(F_{a, K},S^1\cdot x_k^m) \nn\\
&&\qquad = \sum_{1\le k\le q,\ m\ge 1}\dim C_{S^1,d(K)+2N-n-1}(F_{a, K},S^1\cdot x_k^m)\nn\\
&&\qquad = M_{2N-n-1}\ge b_{2N-n-1}=1.\lb{4.52}\eea

Now by (\ref{4.52}) and Lemma 2.2, it yields that there exist at least another closed characteristic $y_n$
with $i(y_n^{2m_n})=2N-n-1$ and $i(y_n^{2m_n})-i(y_n)\in 2\Z$. Thus $y_n$ and its iterations have odd
Maslov-type indices. This completes the proof of Claim 3.

\medskip

Now for Case 2, Theorem 1.2 follows from Claim 2 and Claim 3. The proof of Theorem 1.2 is finished. \hfill\hb

\medskip
{\bf Proof of Theorem 1.3.}

\medskip

Since we are mainly interested in the closed characteristic problem on compact star-shaped hypersurfaces, we just state that our Theorem 1.2 can also be generalized to a class of prequantization bundles based on the work \cite{GGM} of Ginzburg-G\"{u}rel-Macarini, for the reader's convenience, we give some detailed explanations as follows.

The proof is similar to that of Theorem 1.2. It is a generalization of our former work \cite{DLLW} in which we required that every closed characteristic has positive mean index, the proof of Theorem 1.2 is based on our generalization of the common index jump theorem, i.e., Theorem 3.6, and the Morse inequalities (\ref{2.16})-(\ref{2.17}), in which the mean indices of every closed characteristic and its associated symplectic path are only required to be nonzero. In summary, our Theorem 1.2 just weakens the assumption in Theorem 1.2 of \cite{DLLW} that every closed characteristic has positive mean index to that every closed characteristic has nonzero mean index.

While the work \cite{GGM} of Ginzburg-G\"{u}rel-Macarini was also a generalization of our former work \cite{DLLW} to a class of prequantization bundles, the proof of their work also required that every closed orbit has positive mean index, see the first paragraph and the second paragraph in p.408 and p.410 of \cite{GGM} respectively, and also p.438 of \cite{GGM}, then one can notice that their proof is the same as that of \cite{DLLW}.

Now we have improved the common index jump theorem, and note that
there are similar Morse inequalities in the setting of equivariant symplectic homology to
(\ref{2.16})-(\ref{2.17}), cf. (68) in p.221 of \cite{HM}, thus it's natural that we can weaken the condition in Theorems 2.1 and 2.10 of \cite{GGM}
that every closed orbit has positive mean index (or, every closed orbit has negative mean index) to that every closed orbit has nonzero mean index, this is the statement of Theorem 1.3. \hfill\hb

\setcounter{figure}{0}
\setcounter{equation}{0}
\section{Proof of Theorem 1.6} 

In this section, we prove Conjecture 1.5 for the case of $n=3$, i.e., Theorem 1.6, by contradiction. We assume first
the following condition (C):

\medskip

{\bf (C)} {\it There are finitely many prime closed characteristics $\{(\tau_k, y_k)\}_{1\leq k\leq q}$ on $\Sg\in \H_{st}(6)$,
and $\hat{i}(y_k)=0$ for $1\leq k\leq q_0$ with some integers $q_0\in [1,q]$ and $q\in \N$. }

\medskip

Let $P_{\Sg} = \{m\tau_k\;|\;1\le k\le q, m\in\N\}$ be the period set of all closed characteristics on $\Sg$.

Denote by $\ga_k\equiv \ga_{y_k}$ the associated symplectic path of $(\tau_k,y_k)$ for $1\le k\le q$.
Then by Lemma 3.3 of \cite{HuL} and Lemma 3.2 of \cite{Lon1}, there exists $P_k\in Sp(6)$ and $U_k\in Sp(4)$
such that
\bea M_k\equiv\ga_k(\tau_k)=P_k^{-1}(N_1(1,1)\dm U_k)P_k,\qquad \forall\ 1\le k\le q,\lb{5.1}\eea
where, because $\Sigma$ is non-degenerate, every $U_k$ is isotopic in $\Omega^0(U_k)$ to a matrix of the following form by Theorem 3.2 (cf. Theorem 4.7 of \cite{Wan2})
\bea
&& R(\th_1)\,\dm\,\cdots\,\dm\,R(\th_r)\,\dm\,D(\pm 2)^{\dm s} \nn\\
&& \qquad\dm\,N_2(e^{\aa_{1}\sqrt{-1}},A_{1})\,\dm\,\cdots\,\dm\,N_2(e^{\aa_{r_{\ast}}\sqrt{-1}},A_{r_{\ast}})\nn\\
&& \qquad \dm\,N_2(e^{\bb_{1}\sqrt{-1}},B_{1})\,\dm\,\cdots\,\dm\,N_2(e^{\bb_{r_{0}}\sqrt{-1}},B_{r_{0}}), \nn\eea
where $\frac{\th_{j}}{2\pi}\in[0,1]\setminus\Q$ for $1\le j\le r$; $\frac{\aa_{j}}{2\pi}\in[0,1]\setminus\Q$ for $1\le j\le r_{\ast}$;
$\frac{\bb_{j}}{2\pi}\in[0,1]\setminus\Q$ for $1\le j\le r_0$ and
\be r+ s +2r_{\ast} + 2r_0 = 2. \lb{5.2}\ee
Hence by (\ref{5.1}), Theorem 3.8 and the precise index iteration formulae for symplectic paths due to Y. Long
(cf. \cite{Lon3} or Chapter 8 of \cite{Lon4}), for $1\le k\le q$, we have
\bea  i(y_k^m)
&=& m(i(y_k)+n+1-r)+2\sum_{j=1}^r \left[\frac{m\theta_j}{2\pi}\right]+r-1-n  \nn\\
&=& m(i(y_k)+4-r)+2\sum_{j=1}^r \left[\frac{m\theta_j}{2\pi}\right]+r-4,\quad\forall\  m\ge 1, \lb{5.3}\eea
where in (\ref{5.3}), we have used $E(a)=[a]+1$ for $a\in \R\bs\Z$. Thus
\be \hat{i}(y_k)=i(y_k)+4-r+\sum_{j=1}^r\frac{\theta_j}{\pi}, \quad\forall\ 1\le k\le q. \lb{5.4}\ee

\medskip

{\bf Lemma 5.1.} {\it When $\hat{i}(y_k)=0$ with $1\le k\le q_0$, we have $i(y_k^m)=-4$ for any $m\in\N$.}

\medskip

{\bf Proof.} By (\ref{5.4}) we have
\be i(y_k)+4-r+\sum_{j=1}^r\frac{\theta_j}{\pi}=0,\quad \forall\ 1\le k\le q_0. \lb{5.5}\ee
Note that $\frac{\theta_j}{\pi}\notin \Q$ and $0\leq r\leq 2$ by (\ref{5.2}), and then it yields $r=0$ or 2.

If $r=0$, then $i(y_k)=-4$ by (\ref{5.5}). Together with (\ref{5.3}), it yields $i(y_k^m)=-4$ for all $m\in\N$.

If $r=2$, note that $\sum_{j=1}^2\frac{\theta_j}{\pi}\in (0,4)$ and $i(y_k)\in 2\Z$ by (8.1.8) of Theorem 8.1.4
and (8.1.29) of Theorem 8.1.7 of \cite{Lon4}, then by (\ref{5.5}) we have $\sum_{j=1}^2\frac{\theta_j}{\pi}=2$
and $i(y_k)=-4$. Since $\sum_{j=1}^2\frac{m\theta_j}{2\pi}=m$ implies that $\sum_{j=1}^2[\frac{m\theta_j}{2\pi}]=m-1$,
so, by (\ref{5.3}), for any $m\ge 1$, we have
\be i(y_k^m)=-2m+2\sum_{j=1}^r\left[\frac{m\theta_j}{2\pi}\right]-2=-2m+2(m-1)-2=-4. \lb{5.6}\ee

The proof of Lemma 5.1 is finished. \hfill\hb

\medskip

{\bf Remark 5.2.} In our proof of Theorem 1.6 below, we shall apply results in Section 7 of \cite{Vit2} frequently.
Note that in the Theorem 7.1 of \cite{Vit2} and its proof, the two end points of the index interval $I$ were carefully
avoided (cf. Theorem 7.1 on p.637, (7.24) and (7.25) in p.647, and the top part on p.648 of \cite{Vit2}), which are due to
the effect of the $S^1$-action on the homologies with two adjacent dimensions (cf. Corollary in Appendix 1 on p.653 of
\cite{Vit2}). In our proof of Theorem 1.6 we are dealing with only non-degenerate critical orbits. Note that because in our
case the critical orbit $S^1\cdot y_k^m$ studied below is always orientable on the star-shaped hypersurface $\Sg$, by the
homological correspondence of the Theorem I.7.5 and the comments below it on pp.78-79 of \cite{Cha}, only the homology with
dimension to be precisely equal to the Morse index of the critical orbit survives. Therefore in our case, all the results
in (7.24) and (7.25) of \cite{Vit2} work for all $k\in I$, not only for $k\in I^o$. Specially, by (\ref{5.7}) below,
whenever the dimension of the related homology is $d(K')-3$ or $d(K')-5$, i.e., the Viterbo index of the corresponding closed
characteristic is $-3$ or $-5$, results in (7.24) and (7.25) of  \cite{Vit2} can be applied. Such arguments have no
requirement on those Viterbo index near and is not $-3$ and $-5$ of the related closed characteristics, and thus can be
applied to the case in the current paper. This understanding is applied below, whenever we apply Theorem 7.1 of \cite{Vit2}.
We refer readers also to Theorem 5.1 and its proof in \cite{LLW} which generalized Theorem 7.1 of \cite{Vit2} to the
degenerate case.

\medskip

{\bf Proof of Theorem 1.6.}

\medskip

Based on Lemma 5.1, we carry out the proof of Theorem 1.6 in several steps below.

\medskip

{\bf Step 1.} On one hand, by Lemma 5.1, there always holds $i(y_k^m)=-4$ for any $1\le k\le q_0$ and $m\in\N$. On the
other hand, note that $\hat{i}(y_k)>0$ (respectively $<0$) implies $i(y_k^m)\to +\infty$ (respectively $-\infty$) as
$m\to +\infty$. Thus iterates $y_k^m$ of every $y_k$ for $q_0+1\le k\le q$ have indices satisfying $i(y_k^m)\neq -3$ and
$-5$ for any large enough $m\in\N$. Therefore for large enough $a$, all the closed characteristics $y_k^m$ for
$1\le k\le q$ with period larger than $aT$, which implies that the iterate number $m$ is very large, will have their
Viterbo indices:
\bea \left\{\begin{array}{ll}
\mbox{either (i) equal to} -4, \quad \mbox{when}\; \hat{i}(y_k)=0,  \\
\mbox{or (ii) different from} -3,\ -4\ \mbox{and} -5, \quad \mbox{when}\; \hat{i}(y_k)\not= 0.\end{array}\right. \lb{5.7}\eea

\medskip

{\bf Step 2.} For $a\in\R$, let $X^-(a,K)=\{x\in X\mid F_{a,K}(x)<0\}$ with $K=K(a)$ as defined in the above (\ref{2.2})
as well as in Section 7 of \cite{Vit2}. Note also that the origin $0$ of $X$ is not contained in $X^-(a,K)$ by definition.
Because the Hamiltonian function $H_{a,K}$ is quadratic homogeneous as assumed at the beginning of Section 7 of
\cite{Vit2} due to the study there being near the origin, the functional $F_{a,K}$ is homogeneous too.

For any large enough positive $a<a'$, we fix the same constant $K'>0$ as that in the second Step on p.639 of\cite{Vit2}
to be sufficiently large than $K$ such that the Hamiltonian function $H_{t,K'}(x)$ is strictly convex for every
$t\in [a,a']$. Now let $A=X^-(a,K')$ and $A'=X^-(a',K')$.
Because the period set $P_{\Sg}$ defined at the beginning of this section is discrete, we choose the above constants $a$
and $a'$ carefully such that $aT$ and $a'T$ do not belong to $P_{\Sg}$. Note that for $t\in [a,a']$ because every
critical orbit $S^1\cdot x$ of the functional $F_{t,K'}$ always possesses the critical value $F_{t,K'}(S^1\cdot x) = 0$ as
mentioned in p.639 of \cite{Vit2} and by (2.7) of \cite{LLW}, the boundary sets of $A$ and $A'$, i.e.,
$\{x\,|\,F_{t,K'}(x)=0\}$ with $t=a$ or $a'$, possess no critical orbits, and specially the origin $0$ of $X$ is not
contained in $A$ and $A'$. Therefore by the homogeneity mentioned above we have
\bea   H_{S^1, d(K')+i}(A',A) \;=\; H_{S^1, d(K')+i}(A'\cap S(X),A\cap S(X)), \quad \forall\;i\in \Z, \lb{5.8}\eea
where $S(X)$ is the unit sphere of $X$.

Because each critical orbit $S^1\cdot x$ of the functional $F_{t,K'}$ for some $t\in [a,a']$ corresponds to an iterate
$y_k^m$ for some $1\le k\le q$ and $m\in\N$, we denote this critical orbit by $S^1\cdot x_{t,k,m}$. Thus by the definition
of $F_{t,K'}$ we have the period $T_{t,k,m}$ of $x_{t,k,m}$ satisfies $T_{t,k,m}= tT$. Consequently every critical orbit
$S^1\cdot x_{t,k,m}$ of $F_{t,K'}$ for some $t\in [a,a']$ contained in $(A'\bs A)\cap S(X)$ must satisfy
\bea    aT \le T_{t,k,m} \le a'T.   \lb{5.9}\eea
Thus the period of the corresponding $y_k^m$ is also $T_{t,k,m}$ and satisfies (\ref{5.9}) too. Consequently the total
number of such critical orbits contained in $(A'\bs A)\cap S(X)$ is finite. Then we obtain a non-negative integer $\hat{j}$ such that
there exist precisely $\hat{j}$ times of $t\in (a,a')$ which we denote by $t_j$ with $1\le j\le \hat{j}$ satisfying
$a<t_1<\ldots t_{\hat{j}}<a'$, such that $F_{t_j,K'}$ with $1\le j\le \hat{j}$ possesses critical orbits in
$(A'\bs A)\cap S(X)$, and any other $F_{t,K'}$ with $t\in [a,a']\bs \{t_j\,|\,1\le j\le \hat{j}\}$ possesses no
critical orbit in $(A'\bs A)\cap S(X)$.

In order to compute the homology in (\ref{5.8}), we introduce below a new functional $\hat{t}$, which is motivated by
the proof of Proposition 3 in Appendix 1 of \cite{Vit2}.

\medskip

{\bf Claim 4.} {\it The partial derivative $\frac{\pt}{\pt t}F_{t,K'}(x)$ of $F_{t,K'}(x)$ with respect to $t\in [a,a']$
satisfies $\frac{\pt}{\pt t}F_{t,K'}(x)<0$ for all $(x,t)\in S(X)\times [a,a']$.}

\medskip

In fact, by the definition of $H_t(x)$ in Section 2, it is strictly increasing in $t$ when $x\neq0$, and then so is
$H_{t,K'}(x)$. Then the Fenchel dual function $H_{t,K'}^{\ast}(y)$ is strictly decreasing in $t$ when $y\neq0$.
Consequently $F_{t,K'}(x)$ is strictly decreasing in $t$ too. Thus Claim 4 is proved.

Now based on Claim 4 and the well-known implicit function theorem, there exists a unique smooth function
$\hat{t}: S(X)\to \R$ given by the equation
\be F_{\hat{t}(x),K'}(x)=0, \qquad \forall\ x\in S(X).  \lb{5.10}\ee
It further implies
\be \frac{\pt}{\pt t}F_{\hat{t}(x),K'}(x)\hat{t}'(x) + F'_{\hat{t}(x),K'}(x) = 0, \qquad \forall\ x\in S(X). \lb{5.11}\ee
Then for $(x_0,t_0)\in S(X)\times \R$, we have that $x_0$ is a critical point of $\hat{t}$ with critical value $t_0$
if and only if $F'_{t_0,K'}(x_0)=0$ by (\ref{5.11}). Note that $\hat{t}$ is $S^1$-invariant since so is $F_{t,K'}$.

\medskip

{\bf Claim 5.} {\it At any critical point $x_0$ of $\hat{t}$ with critical value $t_0$, we have}
\be C_{S^1, *}(\hat{t}, x_0)\cong C_{S^1, *}(F_{t_0,K'}|_{S(X)}, x_0).\lb{5.12}\ee

\medskip

In fact, let $U$ be a small enough $S^1$-invariant open neighborhood of $S^1\cdot x_0$ in $S(X)$. Since
$\pt F_{t,K'}(x_0)/{\pt t}<0$ by Claim 4, we obtain that $\hat{t}(x)\leq t_0$ for $x\in U$ if and only if $F_{t_0,K'}(x)\leq 0$
by (\ref{5.10}). Thus we obtain
\be \left\{\begin{array}{ll}
\{x\in U\mid\hat{t}(x)\leq t_0 \}=\{x\in U\mid F_{t_0,K'}(x)\leq 0\}, \\
\{x\in U\mid \hat{t}(x)\leq t_0\}\setminus\{S^1\cdot x_0\}=\{x\in U\mid F_{t_0,K'}(x)\leq 0\}\setminus\{S^1\cdot x_0\}. \end{array}\right. \lb{5.13}\ee
Then by the definition of $S^1$-critical module in Sections I.4 and I.7 of \cite{Cha}, the two $S^1$-critical
modules in (\ref{5.12}) are isomorphic to each other at every dimension. Thus Claim 5 holds.

\medskip

{\bf Remark 5.3.} Note that the isomorphic identity (\ref{5.12}) holds without further showing that functional
$\hat{t}(x)$ is $C^2$ and its Morse index and nullity at its critical point $x_0$ with critical value $t_0$ are
the same as those of the functional $F_{t_0,K'}(x_0)$ at its critical point $x_0$, although these can be proved by
using the implicit function theorem and (\ref{5.11}). Here the Hessian matrices of these two functionals differ
by only a positive constant which can be obtained by differentiating both sides of (\ref{5.11}) with respect to $x$,
and then evaluating at the critical points respectively.

\medskip

By Claim 4 and (\ref{5.10}), we then obtain
$$  \left\{\begin{array}{ll}
A'\cap S(X)=\{x\in S(X)\mid F_{a',K'}(x)< 0\}=\{x\in S(X)\mid \hat{t}(x)< a'\}, \\
A\cap S(X)=\{x\in S(X)\mid F_{a,K'}(x)< 0\}=\{x\in S(X)\mid \hat{t}(x)< a\}. \end{array}\right. $$
Note that both $a$ and $a'$ are regular values of $\hat{t}$ since $aT$ and $a'T$ do not belong to $P_{\Sg}$. Then for small enough
$\epsilon>0$, $A'\cap S(X)$ and $A\cap S(X)$ are $S^1$-homotopy equivalent with $\hat{t}^{a'-\epsilon}$ and $\hat{t}^{a+\epsilon}$
respectively, where  $\hat{t}^{\kappa}$ denotes the level set $\hat{t}^{\kappa}=\{x\in S(X)\mid \hat{t}(x)\leq \kappa\}$. Thus by
the homotopy invariance of the homology, we obtain
\be H_{S^1, d(K')+i}(A'\cap S(X),A\cap S(X))\;=\; H_{S^1, d(K')+i}(\hat{t}^{a'-\epsilon},\hat{t}^{a+\epsilon}),
             \forall\;i\in \Z. \lb{5.14}\ee

\medskip

{\bf Step 3.} Now for the chosen large enough $a$ and $a'$ with $a<a'$, by (\ref{5.7}) there exists no any closed
characteristic whose period locates between $aT$ and $a'T$ possessing Viterbo index $-3$ or $-5$. Therefore by the
discussion in pp.78-79 of \cite{Cha}, we obtain
$$  C_{S^1, d(K')-3}(F_{t_j,K'}|_{S(X)}, S^1\cdot x_{t_j,k,m}) = C_{S^1, d(K')-5}(F_{t_j,K'}|_{S(X)}, S^1\cdot x_{t_j,k,m}) = 0, $$
which, together with (\ref{5.12}), yields
\be  C_{S^1, d(K')-3}(\hat{t}, S^1\cdot x_{t_j,k,m}) = C_{S^1, d(K')-5}(\hat{t}, S^1\cdot x_{t_j,k,m}) = 0.   \lb{5.15}\ee
Combining (\ref{5.15}) with an equivariant version of Theorem I.4.3 of \cite{Cha}, i.e., the Morse inequality, we obtain
\be H_{S^1, d(K')-3}(\hat{t}^{a'-\epsilon},\hat{t}^{a+\epsilon})\;=\; H_{S^1, d(K')-5}(\hat{t}^{a'-\epsilon},\hat{t}^{a+\epsilon})\;=\;0. \lb{5.16}\ee
Therefore, combining (\ref{5.16}) with (\ref{5.8}) and (\ref{5.14}), we obtain
\be  H_{S^1, d(K')-3}(A',A) \;=\; H_{S^1, d(K')-5}(A',A) \;=\; 0.   \lb{5.17}\ee

\medskip

{\bf Step 4.} Now we consider the following exact sequence of the triple $(X,A',A)$
\bea
&&\longrightarrow H_{S^1, d(K')-3}(A',A)\stackrel{i_{3*}}{\longrightarrow} H_{S^1, d(K')-3}(X,A)
     \stackrel{j_{3*}}{\longrightarrow}H_{S^1, d(K')-3}(X,A')\nn\\
&&\quad\stackrel{\partial_{3*}}{\longrightarrow} H_{S^1, d(K')-4}(A',A)\stackrel{i_{4*}}{\longrightarrow} H_{S^1, d(K')-4}(X,A)
                         \stackrel{j_{4*}}{\longrightarrow}H_{S^1, d(K')-4}(X,A')\nn\\
&&\qquad\stackrel{\partial_{4*}}{\longrightarrow} H_{S^1, d(K')-5}(A',A)\longrightarrow. \lb{5.18}\eea

Note that by the study of Viterbo in \cite{Vit2}, we have
\bea H_{S^1, d(K)-3}(X,X^-(a,K))
&{\xi_1}\atop{\longrightarrow}& H_{S^1, d(K')-3}(X,X^-(a,K'))  \nn\\
&{\xi_2}\atop{\longrightarrow}& H_{S^1, d(K')-3}(X,X^-(a',K')), \nn\eea
where $\xi_1$ is the homomorphism given by (7.2) of \cite{Vit2}, $\xi_2=j_{3*}$ is the homomorphism given by
the line above (7.4) of \cite{Vit2}. Here the composed homomorphism $\xi=\xi_2\circ\xi_1$ is precisely the
homomorphism given by (7.4) of \cite{Vit2}. But $\xi_1$ is an isomorphism and $\xi$ is a zero homomorphism as
proved in the Steps 1 and 2 of the proof of Theorem 7.1 in \cite{Vit2} respectively. Therefore $j_{3*}=\xi_2$
is also a zero homomorphism.

Therefore (\ref{5.17}) and (\ref{5.18}) yield
\be  H_{S^1, d(K')-3}(X,A) = \mbox{Ker}(j_{3*}) = \mbox{Im}(i_{3*}) = i_{3*}(H_{S^1, d(K')-3}(A',A)) = 0.  \lb{5.19}\ee

Now we fix the above chosen $a'>0$ and choose another large enough $a''>a'$, and enlarge the constant $K'$ in (\ref{2.2})
chosen above (\ref{5.8}) so that the conclusions between (\ref{5.7}) and (\ref{5.8}) hold when we replace $(a,a')$ by
$(a',a'')$. Then repeating the above proof with the long exact sequence of the triple
$(X,A'',A')$ instead of $(X,A',A)$ in the above arguments with $A''=X^-(a'',K')$, similarly we obtain
\be  H_{S^1, d(K')-3}(X,A') =0.  \lb{5.20}\ee

Together with (\ref{5.19}) and (\ref{5.20}), (\ref{5.18}) yields
\be 0\stackrel{\partial_{3*}}{\longrightarrow} H_{S^1, d(K')-4}(A',A)\stackrel{i_{4*}}{\longrightarrow}
   H_{S^1, d(K')-4}(X,A)\stackrel{j_{4*}}{\longrightarrow}H_{S^1, d(K')-4}(X,A')
    \stackrel{\partial_{4*}}{\longrightarrow} 0.  \lb{5.21}\ee

\medskip

{\bf Step 5.} When $a$ increases, we always meet infinitely many closed characteristics with Viterbo index $-4$ due to the
existence of $y_k$ with $\hat{i}(y_k)=0$ for $1\le k\le q_0$ by Lemma 5.1. For the above chosen large enough $a < a'$,
there exist only finitely many closed characteristics among $\{y_k^m\ |\ 1\le k\le q_0, m\ge 1\}$ such that their periods
locate between $aT$ and $a'T$. Therefore for the corresponding critical orbits $S^1\cdot x_{t_j,k,m}$ of $F_{t_j,K'}$, all of
them possess Morse index $d(K')-4$. Then by the equivariant version of Theorem I.4.2 as well as the discussions there on
pp.78-79 of \cite{Cha}, when the condition (C) holds, i.e., $q_0\ge 1$ here, we obtain
\bea  H_{S^1, d(K')-4}(A',A)
&=& H_{S^1, d(K')-4}(A'\cap S(X),A\cap S(X))\nn\\&=& H_{S^1, d(K')-4}(\hat{t}^{a'-\epsilon},\hat{t}^{a+\epsilon})\nn\\
&=&\bigoplus_{aT\le T_{t_j,k,m}\le a'T\atop 1\le k\le q_0}C_{S^1,d(K')-4}(\hat{t},\;S^1\cdot x_{t_j,k,m})  \nn\\
&=&\bigoplus_{aT\le T_{t_j,k,m}\le a'T\atop 1\le k\le q_0}C_{S^1,d(K')-4}(F_{t_j,K'}|_{S(X)},\;S^1\cdot x_{t_j,k,m})  \nn\\
&=& \bigoplus_{aT\le T_{t_j,k,m}\le a'T\atop 1\le k\le q_0}\Q \not= 0, \lb{5.22}\eea
where the first equality follows from (\ref{5.8}), the second equality follows from (\ref{5.14}), the fourth equality follows
from (\ref{5.12}), and for the third equality we give more explanations as follows:

Denote by
\bea
M_q(a,a') &=& \bigoplus_{aT\le T_{t_j,k,m}\le a'T\atop 1\le k\le q}\mathrm{rank}~C_{S^1,q}(\hat{t},\;S^1\cdot x_{t_j,k,m}), \nn\\
\beta_q(a,a') &=& \mathrm{rank}~H_{S^1, q}(\hat{t}^{a'-\epsilon},\hat{t}^{a+\epsilon}). \nn\eea
Then $M_{d(K')-3}(a,a')=M_{d(K')-5}(a,a')=0$, $\beta_{d(K')-3}(a,a')=\beta_{d(K')-5}(a,a')=0$ by (\ref{5.15}) and
(\ref{5.16}) respectively, which together with an equivariant version of Theorem I.4.3 of \cite{Cha} yield
$M_{d(K')-4}(a,a')=\beta_{d(K')-4}(a,a')$. Then the third equality in (\ref{5.22}) holds.

\medskip

{\bf Step 6.} By the exactness of the sequence (\ref{5.21}) and (\ref{5.22}), we obtain
$$   H_{S^1, d(K')-4}(X,A)=H_{S^1, d(K')-4}(A',A)\bigoplus H_{S^1, d(K')-4}(X,A')\neq 0.  $$
Then, by our choice of $a$, $a'$ and $a''$, and replacing $(X,A',A)$ by $(X,A'',A')$ in the above arguments, similarly
we obtain
\be H_{S^1, d(K')-4}(X,A')\neq 0. \lb{5.23}\ee

Now on one hand, if $j_{4*}$ in (\ref{5.21}) is a trivial homomorphism, then by the exactness of the sequence
(\ref{5.21}) it yields
$$   H_{S^1, d(K')-4}(X,A') = \mbox{Ker}(\partial_{4*}) = \mbox{Im}(j_{4*})=0,   $$
which contradicts to (\ref{5.23}). Therefore $j_{4*}$ in (\ref{5.21}) is a non-trivial homomorphism.

However, on the other hand, by the same proof of the fact that $j_{3*}$ is a zero homomorphism in Step 4, we obtain $j_{4*}$ in (\ref{5.21}) is also a zero homomorphism. This contradiction completes the proof of Theorem 1.6. \hfill\hb

\vspace{5mm}

\noindent {\bf Acknowledgment}

\vspace{2mm}

The authors would like to express their great gratitude to the referee for her/his careful reading, many
valuable comments and suggestions, which greatly improved this manuscript.

\bibliographystyle{abbrv}

\end{document}